\input amstex.tex
\input amsppt.sty   
\magnification 1200
\vsize = 9.5 true in
\hsize=6.2 true in
\NoRunningHeads        
\parskip=\medskipamount
        \lineskip=2pt\baselineskip=18pt\lineskiplimit=0pt
       
        \TagsOnRight
        \NoBlackBoxes

        \topmatter
        \title
        Pure Point Spectrum of the Floquet Hamiltonian \\for 
        the Quantum Harmonic Oscillator \\Under Time Quasi- 
        Periodic Perturbations       
        \endtitle
        \author
        W.-M.~Wang
        \footnote"{$^{*}$}"{Partially supported by NSF grant
DMS-05-03563\hfill\hfill\qquad}
        \endauthor
\address
Departement de Mathematique, Universite Paris Sud, 91405 Orsay Cedex, FRANCE;
Department of Mathematics, University of Massachussetts, Amherst, Ma 01003, USA
\endaddress
        \email
{wei-min.wang\@math.u-psud.fr,  weimin\@math.umass.edu}
\endemail
\thanks  
{I am deeply indebted to the referee for a thorough reading of the paper and for many helpful comments.}
\endthanks 
\abstract
We prove that the $1-d$ quantum harmonic oscillator is stable under spatially localized,
time quasi-periodic perturbations on a set of Diophantine frequencies of positive measure.
This proves a conjecture raised by Enss-Veselic in their 1983 paper \cite{EV} in the general
quasi-periodic setting.  The motivation of the present paper also comes from construction
of quasi-periodic solutions for the corresponding nonlinear equation.  
\endabstract

        \bigskip\bigskip
        \bigskip
        \toc
        \bigskip
        \bigskip 
        \widestnumber\head {Table of Contents}
        \head 1. Introduction and statement of the theorem\endhead
        \head 2. The Floquet Hamiltonian in the Hermite-Fourier  basis 
        \endhead
        \head 3. Exponential decay of Green's functions at fixed $E$: estimates in $\theta$
        \endhead
        \head 4. Frequency estimates and the elimination of $E$
        \endhead
        \head 5. Proof of the Theorem
        \endhead
        \endtoc
        \endtopmatter
        \vfill\eject
        \bigskip
\document
\head{\bf 1. Introduction and statement of the theorem}\endhead
The stability of the quantum harmonic oscillator is a long standing problem since
the establishment of quantum mechanics. The Schr\"odinger equation for the harmonic oscillator
in $\Bbb R^n$ (in appropriate coordinates) is the following:
$$-i\frac\partial{\partial t}\psi =\frac{1}{2}\sum_{i=1}^{n}(-\frac{\partial^2}{\partial x_i^2}+x_i^2)\psi,
\tag 1.1$$
where we assume 
$$\psi\in C^1(\Bbb R, L^2(\Bbb R^n))$$
for the moment. We start from the $1$ dimensional case, $n=1$. (1.1) then reduces to
$$-i\frac\partial{\partial t}\psi =\frac{1}{2}(-\frac{\partial^2}{\partial x^2}+x^2)\psi.
\tag 1.2$$
The Schr\"odinger operator 
$$H{\overset\text{def }\to =}(-\frac{d^2}{d x^2}+x^2)\tag 1.3$$
is the 1-d harmonic oscillator. Since $H$ is independent of $t$, it is amenable to a spectral
analysis. It is well known that $H$ has pure point spectrum with eigenvalues
$$\lambda_n=2n+1, \quad n=0,\,1..,\tag 1.4$$
and eigenfunctions (the Hermite functions)
$$h_n(x)=\frac{H_n(x)}{\sqrt {2^n n!}}e^{-x^2/2},\quad n=0,\,1...\tag 1.5$$
where $H_n(x)$ is the $n^{\text {th}}$ Hermite polynomial, relative to the weight
$e^{-x^2}$ ($H_0(x)=1$) and
$$\aligned &\int_{-\infty}^{\infty}e^{-x^2}H_m(x)H_n(x)dx\\
=&2^n n!\sqrt\pi\delta_{mn}.\endaligned\tag 1.6$$
Using (1.4-1.6), the normalized $L^2$ solutions to (1.1) are all of the form
$$\psi(x,t)=\sum_{n=0}^{\infty}a_nh_n(x)e^{i\frac{\lambda_n}{2}t}\quad (\sum|a_n|^2=1),\tag 1.7$$
corresponding to the initial condition $$\psi(x,0)=\sum_{n=0}^{\infty}a_nh_n(x)\quad(\sum|a_n|^2=1).\tag 1.8$$
The functions in (1.7) are almost-periodic (in fact periodic here) in time with frequencies $\lambda_n/4\pi,\,n=0,\,1...$

Equation (1.2) generates a unitary propagator $U(t, s)=U(t-s, 0)$ on $L^2(\Bbb R)$. Since the spectrum of $H$
is pure point, $\forall u\in L^2(\Bbb R)$, $\forall\epsilon$, $\exists R$, such that
$$\inf_t\Vert U(t,0)u\Vert_{L^2(|x|\leq R)}\geq (1-\epsilon)\Vert u\Vert\tag 1.9$$
by using eigenfunction (Hermite function) expansions.
\smallskip
The harmonic oscillator (1.3) is an integrable system. The above results are classical. It is natural to
ask how much of the above picture remains under perturbation, when the system is no longer integrable.
In this paper, we investigate stability of the 1-d harmonic oscillator under time quasi-periodic, spatially
localized perturbations. To simplify the exposition, we study the following ``model" equation:
$$-i\frac\partial{\partial t}\psi
=\frac{1}{2}(-\frac{\partial^2}{\partial x^2}+x^2)\psi+\delta|h_0(x)|^2\sum_{k=1}^\nu\cos(\omega_kt+\phi_k)\psi,
\tag 1.10$$
on $C^1(\Bbb R, L^2(\Bbb R))$, where
$$\aligned &0<\delta\ll 1,\\
&\omega=\{\omega_k\}_{k=1}^\nu\in[0, 2\pi)^\nu,\\
&\phi=\{\phi_k\}_{k=1}^\nu\in[0, 2\pi)^\nu,\\
&h_0(x)=e^{-x^2/2}\quad \text{is the 0th Hermite function.}\endaligned\tag 1.11$$
In particular, we shall study the validity of (1.9) for solutions to (1.10), when $U$ is the propagator 
for (1.10). The method used here can be generalized to treat the equation
$$ -i\frac\partial{\partial t}\psi
=\frac{1}{2}(-\frac{\partial^2}{\partial x^2}+x^2)\psi+\delta V(t,x),$$ 
where $V$ is $C_0^\infty$ in $x$ and analytic, quasi-periodic in $t$. 

The perturbation term, $\Cal O(\delta)$ term in (1.10) is motivated by the nonlinear equation:
$$-i\frac\partial{\partial t}\psi
=\frac{1}{2}(-\frac{\partial^2}{\partial x^2}+x^2)\psi+M\psi
+\delta|\psi|^2\psi\quad (0<\delta\ll 1),
\tag 1.12$$
where $M$ is a Hermite multiplier, i.e., in the Hermite function basis, 
$$\aligned M&=\text {diag}\,(M_n),\quad M_n\in\Bbb R,\\
Mu&=\sum_{n=0}^\infty M_n(h_n, u)h_n,\text{   for all } u\in L^2(\Bbb R).\endaligned$$
Specifically, (1.10) is motivated by the construction of time quasi-periodic
solutions to (1.12) for appropriate initial conditions such as 
$$\psi (x,0)=\sum_{i=1}^\nu c_{k_i}h_{k_i}(x).\tag 1.13$$ In (1.10), for computational simplicity,
we take the spatial dependence to be $|h_0(x)|^2$ as it already captures the essence of the perturbation
in view of (1.12, 1.13, 1.5). The various computations and the Theorem extend immediately to more general finite combinations of $h_k(x)$.

\noindent{\it The Floquet Hamiltonian and formulation of stability}

It follows from \cite{Y2, 3} that (1.10) generates a unique unitary propagator $U(t,s)$, $t,\,s\in\Bbb R$
on $L^2(\Bbb R)$, so that for every $s\in\Bbb R$ and $$u_0\in B^2=\{f\in L^2(\Bbb R)|
\Vert f\Vert^2_{B^2}=\sum_{|\alpha+\beta|\leq 2}\Vert x^\alpha\partial_x^\beta f\Vert ^2_{L^2}<\infty\},$$
$$u(\cdot)=U(\cdot, s)u_0\in C^1(\Bbb R,\,L^2(\Bbb R))\cap C^0(\Bbb R,\,B^2)\tag 1.14$$
is a unique solution of (1.10) in $L^2(\Bbb R)$ satisfying $u(s)=u_0$.

When $\nu=1$, (1.10) is time periodic with period $T=2\pi/\omega$. The 1-period propagator
$U(T+s, s)$ is called the Floquet operator. The long time behavior of the solutions
to (1.10) can be characterized by means of the spectral properties of $U(T+s, s)$ \cite{EV,Ho,YK}.
Furthermore the nature of the spectrum of $U$ is the same (apart from multiplicity) as that of the
Floquet Hamiltonian $K$ \cite{Y1}:
$$K=i\omega \frac{\partial}{\partial\phi}+\frac{1}{2}(-\frac{\partial^2}{\partial x^2}+x^2)\psi+\delta|h_0(x)|^2\cos\phi$$
on $L^2(\Bbb R)\otimes L^2(\Bbb T)$, where $L^2(\Bbb T)$ is $L^2[0, 2\pi)$ with periodic boundary 
conditions. 

Decompose $L^2(\Bbb R)$ into the pure point $\Cal H_{pp}$ and continuous $\Cal H_c$
spectral subspaces of the Floquet operator $U(T+s,s)$:
$$L^2(\Bbb R)=\Cal H_{pp}\oplus\Cal H_c.$$
We have the following equivalence relations \cite{EV,YK}: 
$u\in\Cal H_{pp}(U(T+s,s))$ if and only if
$\forall \epsilon>0$, $\exists R>0$, such that 
$${\inf}_t\Vert U(t,s)u\Vert_{L^2(|x|\leq R)}\geq (1-\epsilon)\Vert u\Vert;$$
and 
$u\in\Cal H_{c}(U(T+s,s))$ if and only if 
$\forall R>0$,  
$${\lim}_{t\to\pm\infty}\frac{1}{t}\int_0^tdt'\Vert U(t',s)u\Vert^2_{L^2(|x|\leq R)}=0.$$
(Needless to say, the above statements hold for general time periodic Schr\"odinger equations.)

 When $\nu\geq 2$, (1.10) is time quasi-periodic. The above constructions extend for small $\delta$,
 cf. \cite{Be, E, JL} leading to the Floquet Hamiltonian $K$:
 $$K=i\sum_{k=1}^\nu\omega_k\frac{\partial}{\partial\phi_k}+\frac{1}{2}(-\frac{\partial^2}{\partial x^2}+x^2)\psi+\delta|h_0(x)|^2\sum_{k=1}^\nu\cos\phi_k\tag 1.15$$
 on $L^2(\Bbb R)\otimes L^2(\Bbb T^\nu)$, cf. \cite{BW1}. 
 This is related to the so called reducibility of skew product flows in dynamical systems, cf. \cite {E}.
 We note that the Hermite-Fourier functions:
 $$e^{-in\cdot\phi}h_j(x), \quad n\in\Bbb Z^\nu,\quad \phi\in\Bbb T^\nu,\quad j\in\{0,\,1...\}\tag 1.16$$ provide a basis for $L^2(\Bbb R)\otimes L^2(\Bbb T^\nu)$.

We say that the harmonic oscillator $H$ is {\it stable} if $K$ has pure point spectrum. Let
$s\in\Bbb R$.  This implies (by
expansion using eigenfunctions of $K$) that
given any $u\in L^2(\Bbb R)$, $\forall \epsilon>0$, $\exists R>0$, such that 
$${\inf}_t\Vert U(t,s)u\Vert_{L^2(|x|\leq R)}\geq (1-\epsilon)\Vert u\Vert, \text{a.e.}\,\phi,\tag 1.17$$
cf. \cite{BW1, JL}. So (1.9) remains valid and we have dynamical stability.
\smallskip
We now state the main results pertaining to (1.10).
\proclaim {Theorem}
There exists $\delta_0>0$, such that for all $0<\delta<\delta_0$, there exists $\Omega\subset [0,2\pi)^\nu$
of positive measure, asymptotically full measure:
$$\text{mes }\Omega\to(2\pi)^\nu\quad \text{as }\delta\to 0,$$
such that 
for all $\omega\in\Omega$, the Floquet Hamiltonian $K$ defined in (1.15) has pure point spectrum: 
$\sigma(K)=\sigma_{\text {pp}}$. Moreover the Fourier-Hermite coefficients of the eigenfunctions 
of $K$ have subexponential decay.
\endproclaim

As an immediate consequence, we have
\proclaim{Corollary}
Assume that $\Omega$ is as in the Theorem. Let $s\in\Bbb R$. For all $\omega\in\Omega$, all $u\in L^2(\Bbb R)$, all $\epsilon>0$,
there exists $R>0$, such that 
$${\inf}_t\Vert U(t,s,\phi)u\Vert_{L^2(|x|\leq R)}\geq (1-\epsilon)\Vert u\Vert, \text{a.e.}\,\phi,\tag 1.18$$
where $U$ is the unitary propagator for (1.10).
\endproclaim

We note that this good set  $\Omega$ of $\omega$ is a subset of Diophantine frequencies. This is typical for KAM type of persistence theorem. Stability under time quasi-periodic perturbations as in (1.10) is,
generally speaking a precursor for stability under nonlinear perturbation as in (1.12) (cf. \cite{BW1, 2}),
where $M$ plays the role of $\omega$ and varies the tangential frequencies.  The above Theorem
resolves the Enss-Veselic conjecture dated from their 1983 paper \cite{EV} in a general quasi-periodic
setting.
\smallskip
\noindent{\it A sketch of the proof of the Theorem.}

Instead of working with $K$ defined on $L^2(\Bbb R)\otimes L^2(\Bbb T^\nu)$
directly, it is more convenient to work with its unitary equivalent $H$ on $\ell^2(\Bbb Z^\nu\times\{0,\,1...\})$,
using the Hermite-Fourier basis in (1.16). We have 
$$H=\text {diag }(n\cdot\omega+j+\frac{1}{2})+\frac{\delta}{2}W\otimes\Delta\tag 1.19$$
on $\ell^2(\Bbb Z^\nu\times\{0,\,1...\})$, where 
$W$ acts on the $j$ indices, $j=0,\,1,\,2...$,
$$W_{jj'}\sim\frac{1}{\sqrt{j+j'}}e^{-\frac{(j-j')^2}{2(j+j')}}\quad \text {for } j+j'\gg 1;\tag 1.20$$
$\Delta$ acts on the $n$ indices, $n\in\Bbb Z^\nu$,
$$\aligned \Delta_{nn'}&=1,\quad |n-n'|_{\ell^1}=1,\\
&=0,\quad \text{otherwise}.\endaligned\tag 1.21$$

The reduction from (1.15) to (1.19-1.21) is performed in section 2. The main work is to compute 
$W$, which involves integrals of products of Hermite functions. We will explain shortly
this computation, which is independent from the main thread of construction.

The principal new feature here is that $W$ is long range. The $j^{\text{th}}$ row has width 
$\Cal O(\sqrt j)$ about the diagonal element $W_{jj}$. It is  
{\it not} and {\it cannot}  be approximated by a convolution matrix. The potential $x^2$ breaks translational invariance. The annihilation and creation 
operators of the harmonic oscillator
$a=\frac{1}{\sqrt 2}(\frac{d}{dx}+x)$,  $a^*=\frac{1}{\sqrt 2}(-\frac{d}{dx}+x)$, satisfying $[a, a^*]=1$, are
generators of the Heisenberg group. So (1.19) presents a new class of problems distinct from that
considered in \cite{B1-3, BW1, 2, EK, Ku1, KP}. 

The proof of pure point spectrum of $H$ is via proving pointwise decay as $|x-y|\to\infty$ 
of the finite volume Green's functions: $(H_\Lambda-E)^{-1}(x, y)$, 
where $\Lambda$ are finite subsets of $\Bbb Z^\nu\times\{0,\,1...\}$
and $\Lambda\nearrow \Bbb Z^\nu\times\{0,\,1...\}$. We need decay of the Green's functions at
all scales, as assuming $E$ an eigenvalue, a priori we do not have information on the center
and support of its eigenfunction $\psi$. The regions $\Lambda$ where $(H_\Lambda-E)^{-1}$
has pointwise decay is precisely where we establish later that $\psi$ is small there.

For the initial scales, the estimates on $G_\Lambda(E)=(H_\Lambda-E)^{-1}$ are obtained
by direct perturbation theory in $\delta$ for $0<\delta\ll 1$. For subsequent scales, the proof
is a multiscale induction process using the resolvent equation. Assume we have estimates on 
$G_{\Lambda'}$ for cubes $\Lambda'$ at scale $L'$ and $\Lambda$ is a cube at a larger
scale $L$, $L\gg L'$. Intuitively, if we could establish that for most of $\Lambda'\subset\Lambda$,
$G_{\Lambda'}(E)$ has pointwise decay, then assuming we have some a priori estimates on 
$G_\Lambda(E)$, we should be able to prove that $G_\Lambda(E)$ also has pointwise decay. 

There are ``two" directions in the problem, the higher harmonics direction $n$ and the spatial 
direction $j$. The off-diagonal part of $H$ is Toeplitz in the $n$ direction, corresponding to the 
discrete Laplacian $\Delta$. Since the frequency $\omega$ is in general a vector (if $\nu\geq 2$),
$n\cdot\omega$ does not necessarily $\to\infty$ as $|n|\to\infty$. So the $n$ direction is non 
perturbative. We use estimates on $G_{\Lambda'}$ and semi-algebraic techniques as in 
\cite{BGS,÷BW1} to control
the number of resonant $\Lambda'$, where $G_{\Lambda'}$ is large, in $\Lambda$. 

In the $j$ direction, we do analysis, i.e., perturbation theory. This is the new feature.
From (1.19) and Schur's lemma, $\Vert W\otimes\Delta\Vert=\Cal O(1)$. So the $\ell^2$ norm 
of the perturbation does not  decay (relative to eigenvalue spacing) in $j$. However when $\delta=0$, $H$ is diagonal with eigenvalues $n\cdot\omega+j$ and eigenfunctions $\delta_{n,j}$, the canonical basis for $\ell^2(\Bbb Z^\nu\times\{0,\,1...\})$. We have 
$$\Vert [W\otimes\Delta]\delta_{n,j}\Vert =\Cal O(\frac{1}{j^{1/4}})\quad (j\geq 1),$$
which decays in $j$. 

This is intuitively reasonable, as $W$ stems from a spatially localized perturbation from (1.10).
As $j$ increases, The Hermite functions $h_j$ become more extended, c.f. (1.5). So the effect
of the spatial perturbation should decrease as $j$ increases. 

Assuming $\omega$ is Diophantine: $$\Vert n\cdot\omega\Vert_{\Bbb T}\geq\frac{c}{|n|^\alpha}\quad (c>0,\,n\neq 0,\,\alpha>2\nu),$$
where $\Vert\cdot\Vert_{\Bbb T}$ is the distance to the nearest integer,
this enables us to preserve local eigenvalue spacing for $\Lambda$
which are appropriately proportioned in $n$, $j$. This in turn leads to decay of Green's functions.
Combining the estimates in the $n$ and $j$ directions, we obtain estimates on the Green's
function at the larger scale $L$.

\smallskip
\noindent{\it Integrals of products of Hermite functions.}

From (1.15, 1.19), computation of $W$ involves computing the following integrals:
$$\aligned &\int_{-\infty}^{\infty}h_0^2(x)h_m(x)h_n(x)dx\\
=&\frac{1}{\sqrt{2^{n+m}m!n!}}\int_{-\infty}^{\infty}e^{-2x^2}H_0^2(x)H_m(x)H_n(x)dx,\quad m,\,n=0,\,1...,\endaligned\tag 1.22$$
where $H_m$, $H_n$ are respectively the $m^{\text {th}}$, $n^{\text {th}}$ Hermite polynomial, 
$H_0(x)=1$.

Let $$I=\int_{-\infty}^{\infty}e^{-2x^2}H_0^2(x)H_m(x)H_n(x)dx.$$ The idea is to view $e^{-x^2}H_0^2(x)$
as $ e^{-x^2}H_0(\sqrt 2 x)$, i.e., the $0^{\text {th}}$ Hermite function relative to the weight $e^{-2x^2}$
and to use the generating function of Hermite polynomials to reexpress 
$$H_m(x)H_n(x)=\sum_{\ell=0}^{m+n}a_\ell H_\ell(\sqrt 2 x).\tag 1.23$$
We then have $$\aligned I&=a_0\int[H_0(\sqrt 2 x)]^2e^{-2x^2}dx\\
&=a_0\sqrt{\pi/2}\endaligned$$
using (1.6). This computation is carried out in section 2, (2.7-2.10), recovering an apparently classical
result, which could be found in e.g., \cite{GR, PBM}. 
\smallskip
More generally, we are interested in computing 
$$I=\int_{-\infty}^{\infty}e^{-2x^2}H_p(x)H_q(x)H_m(x)H_n(x)dx,\quad p,\,q,\,m,\,n=0,\,1... \tag 1.24$$
which are needed for the nonlinear equation or if we consider more general perturbations of the
harmonic oscillator.
Following the same line of arguments, we decompose $H_p(x)H_q(x)$ into  
$$H_p(x)H_q(x)=\sum_{\ell=0}^{p+q}b_\ell H_\ell(\sqrt 2 x).\tag 1.25$$
Combining (1.23) with (1.25), and assuming (without loss of generality), $p+q\leq m+n$, we then have
$$I=\sum_{\ell=0}^{p+q}a_\ell b_\ell c_\ell,$$
where $c_\ell=\int_{-\infty}^{\infty}[H_\ell(\sqrt 2 x)]^2e^{-2x^2}dx.$

The computation for general $p$, $q$ is technically more involved and is carried out in 
\cite{W}. Unlike the special case $p=q=0$, we did not find the corresponding result for
general $p$, $q$ in existing literature. 

The computation of $I$ in (1.24) is {\it exact} (see (2.10)), reflecting the integrable nature
of the quantum harmonic oscillator. The proof of the theorem is, however, general. It is
applicable as soon as the kernel $W$ satisfies (1.20). 
Following the precedent discussion on $I$ for general $p$, $q$, and using properties of the Hermite series 
(cf. \cite{T} and references therein), one should be able to extend the Theorem to $V$, which are 
 $C_0^\infty$in $x$ and analytic quasi-periodic in $t$, leading to perturbation kernels in the
 Hermite-Fourier basis satisfying conditions similar to (1.20) in the $j$ direction  and exponential
 decay condition in the $n$ direction. 

When the perturbation $V$ is independent of time and is a $0^{\text{th}}$ order symbol, satisfying
$$|\partial^\alpha V|\leq C_\alpha (1+|x|)^{-\alpha},\quad \alpha=0,\,1...\tag 1.26$$
the corresponding Schr\"odinger equation has been studied in e.g., \cite{BBL, KRY, Z},
where it was shown that certain properties of the harmonic oscillator equation extend to the perturbed 
equation. The spectral property needed for the construction here is more detailed and
stringent. Hence it is reasonable to believe that the set of potentials $V$ will be more restrictive than that
in (1.26)

\noindent {\it Some perspectives on the Theorem}

The Theorem shows that for small $\delta$, there is a subset $\Omega\subset[0,2\pi)^\nu$ of Diophantine frequencies of
positive measure, such that if $\omega\in\Omega$, then (1.9) is satisfied.
Hence spatially localized solutions remain localized for all time. It is natural to ask what happens if the forcing
frequencies $\omega$ are in the complement set, $\omega\in\Omega^c$.

If $\omega$ is rational, the perturbing potential $V$ is bounded and has sufficiently fast 
decay at infinity, it is known from general compactness argument \cite{EV} that the Floquet
Hamiltonian has pure point spectrum. In our example, this can be seen as follows. 
In (1.10) restricting to periodic perturbation ($\nu=1$), it is easy to see that for
$$\forall \omega,\,A=(n\omega+j+z)^{-1} W\otimes\Delta,\text{ where } \Im z=1\text { is compact}.\tag 1.27 $$

Assume $\omega$ is rational: $\omega=p/q,\, (q\neq 0)$. Since $H_0=n\omega+j$ has pure point
spectrum (with infinite degeneracy) and the spacing between different eigenvalues is $1/q$, 
(1.27) implies that $H$ has pure point spectrum. When $\omega$ is irrational,
$H_0$ typically has dense spectrum. No conclusion can be drawn from (1.27).  

It is worth remarking that (1.27) holds {\it for all scalar} $\omega$. In the quasi-periodic case, $\omega$
is a vector, the compactness argument breaks down. The proof of Lemma 3.5 in the present paper
is a replacement.

If $V$ is unbounded, we have a different situation. The results in \cite{HLS, GY} combined with \cite{YK} show that
for the following {\it unbouded} time {\it periodically} perturbed harmonic oscillator:
$$i\frac{\partial u}{\partial t}=\frac{1}{2}(-\Delta+x^2)u+2\epsilon (\sin t)x_1 u+\mu V(t,x) u,\quad
x=(x_1,... x_n)\in\Bbb R^n\tag 1.28$$
where $V(t,x)$ is a real valued smooth function of $(t,x)$, satisfying
$$V(t+2\pi, x)=V(t),\quad |V(t,x)|\lesssim|x|\quad\text {as }x\to\infty,\quad |\partial_x^{\alpha}V(t,x)|\leq
C_\alpha,\quad |\alpha|\geq 1,$$ 
the solutions diffuse to infinity as $t\to\infty$. More precisely, for all 
$u_0\in L^2(\Bbb R^n)\cap H^2 (\Bbb R^n)$, for any $R>0$, the solution $u_t$
satisfies
$$\lim_{T\to\pm\infty}\frac{1}{T}\int_0^Tdt\Vert u_t\Vert_{L^2(|x|\leq R)}=0.\tag 1.29$$

In (1.28), $\nu=1$ (periodic), $\omega=1$ , $\omega\in\Omega^c$, (1.29) is an opposite of (1.9).
However the perturbation is unbounded. Moreover the proof in \cite{GY} uses in an essential way that the potential is linear at infinity,  hence positivity of the commutator: $[\frac{d}{dx_1}, x_1]=1$.  

In the exactly solvable case where the time periodic perturbations is quadratic in the spatial
coordinates, it is known that the Floquet  Hamiltonian exhibits a transition between pure
point and continuous spectrum as the frequency is varied \cite{Co1}. The perturbation there
is again unbounded. 
\smallskip
\noindent {\it Some related results.}

To our knowledge, when $\omega\in\Omega$ is nonresonant, there were no results in the literature on the perturbed harmonic oscillator equation of type (1.10), even in the time periodic case, i.e., $\omega\in
[0,2\pi)$. The main difficulties encountered by the traditional KAM method seem to be (i) the eigenvalue
spacing for the unperturbed operator does {\it not} grow, $\lambda_{k+1}-\lambda_k=1$, (ii) the perturbation $W$ in the Hermite basis has slow decay (1.20).

When the eigenvalue spacing for the unperturbed operator grows: $|E_{j+1}-E_j|>j^\beta\quad (\beta>0)$, which
corresponds to a potential growing faster than quadratically at infinity, and when the perturbation
is {\it periodic} in time, related stability results were proven in \cite{DS}. In \cite{Co2}, 
under time {\it periodic} perturbation and replacing $W$ in (1.20) by a faster decaying kernel, hence 
decaying norm in $j$, which no longer corresponds to the physical case of harmonic oscillator under time periodic, spatially localized perturbation, stability results were also proven. Both papers used some modified
KAM method.    
\smallskip 
\noindent {\it Motivation for studying (1.10)} 

As mentioned earlier, the motivation for analyzing (1.10) partly comes from the nonlinear equation (1.12).
In \cite{B1-3, EK}, time quasi-periodic solutions were constructed for the nonlinear Schr\"odinger equation
in $\Bbb R^d$ with Dirichlet or periodic boundary condition
$$i\frac\partial{\partial t}\psi =(-\Delta+M)\psi+\delta |\psi|^{2p}\psi,\qquad (p\in\Bbb N^+;\,0<\delta\ll
1)\tag 1.30
$$
where $M$ is a Fourier multiplier; see \cite{Ku1, KP} for the Dirichlet case in $\Bbb R$ with a potential
in place of $M$.  In
\cite{BW2}, time quasi-periodic solutions were constructed for the nonlinear random Schr\"odinger equation
in $\Bbb Z^d$
$$i\frac\partial{\partial t}\psi =(-\epsilon\Delta+V)\psi+\delta |\psi|^{2p}u,\qquad (p\in\Bbb
N^+\,;\epsilon,\,0<\delta\ll 1),\tag 1.31 
$$
where $V=\{v_j\}_{j\in\Bbb Z^d}$ is a family of random variables. 

The proofs in \cite{B1-3, BW2} use 
operator method, which traces its origin to the study of Anderson localization \cite{FS}. This method
was first applied in the context of Hamiltonian PDE in \cite{CW}. The proofs in \cite{EK, Ku1, KP}
use KAM type of method.  

In (1.30) (specializing to 1-d), the eigenvalues of the linear operator are $n^2$, so 
$E_{n+1}-E_n\sim n$, the eigenfunctions $e^{inx}$, however, are extended: $|e^{inx}|=1$ for all $x$. 
Let us call this case $A$, where there is eigenvalue separation.
In (1.31), the eigenvalues of the linear operator form a {\it dense} set, the eigenfunctions, on the other
hand are not only {\it localized} but localized about {\it different} points in $\Bbb Z^d$ from Anderson
localization theory, see e.g., \cite{GB, GK}. This is case $B$, where there is eigenfunction separation.
The existence of time quasi-periodic solutions, i.e., KAM type of solutions in $A$ is a consequence
of eigenvalue separation; while in $B$, eigenfunction separation.

Equation (1.10) and its nonlinear counterpart
$$-i\frac\partial{\partial t}\psi =\frac{1}{2}(-\frac {\partial^2}{\partial x^2}+x^2)\psi+M\psi+\delta
|\psi|^{2p}\psi,\qquad (p\in\Bbb N^+;\,0<\delta\ll 1),\tag 1.32
$$
where $M$ is a Hermite multiplier, stand apart from both (1.30, 1.31). It is neither $A$, nor $B$.
There is  eigenvalue spacing, but it is a constant: $\lambda_{n+1}-\lambda_n=1$. In particular, it does not grow
with $n$. The eigenfunctions (Hermite functions) $h_n$ are  ``localized" about the origin. But they
become more extended as $n$ increases because of the presence of the Hermite polynomials, cf. (1.5). 
This in turn leads to the long range kernel $W$ in (1.20) and long range nonlinearity in (1.32) in the 
Hermite function basis, cf. \cite{W}. 

From the KAM perspective a la Kuksin, this is a borderline case, where Theorem 1.1
in \cite {Ku2} does not apply. The more recent KAM type of theorem in \cite{EK} does not apply 
either, because $W$ is long range and not close to a Toeplitz matrix (cf.(1.20)) for the reasons stated earlier. These are the features which make (1.10, 1.32) interesting from a mathematics point of view, aside from its apparent relevance to physics.
\bigskip
\head{\bf 2. The Floquet Hamiltonian in the Hermite-Fourier basis}\endhead
Recall from section 1, the Floquet Hamiltonian
$$K={i}\sum_{k=1}^\nu\omega_k\frac{\partial}{\partial\phi_k}+\frac{1}{2}(-\frac{\partial^2}{\partial x^2}+x^2)+\delta|h_0(x)|^2\sum_{k=1}^\nu\cos\phi_k\tag 2.1$$
on $L^2(\Bbb R)\otimes L^2(\Bbb T^\nu)$, where
$$\aligned &0<\delta\ll 1,\\
&\omega_k\in[0, 2\pi),\quad k=1,...,\nu\\
&\phi_k\in[0, 2\pi),\quad k=1,...,\nu\\
&h_0(x)=e^{-x^2/2}.\endaligned\tag 2.2$$
As mentioned in section 1, $h_0(x)$
is the $0^{\text {th}}$ Hermite function, $0^{\text {th}}$ eigenfunction of the 1-d harmonic oscillator
and more generally,
$$\aligned (-\frac{d^2}{d x^2}+x^2)h_n&=\lambda_nh_n,\\
\lambda_n&=2n+1,\quad n=0,\,1...\\
h_n(x)&=\frac{H_n(x)}{\sqrt {2^n n!}}e^{-x^2/2},\quad n=0,\,1...,\endaligned\tag 2.3$$
where $H_n(x)$ is the $n^{\text {th}}$ Hermite polynomial, relative to the weight
$e^{-x^2}$ ($H_0(x)=1$) and
$$\aligned &\int_{-\infty}^{\infty}e^{-x^2}H_m(x)H_n(x)dx\\
=&2^n n!\sqrt\pi\delta_{mn}\endaligned\tag 2.4$$

\noindent{\it Integral of products of Hermite functions}

We express (2.1) in the Hermite function basis for small $\delta$ and 
compute the integral:
$$\aligned &\int_{-\infty}^{\infty}h_0^2(x)h_m(x)h_n(x)dx\\
=&\frac{1}{\sqrt{2^{n+m}m!n!}}\int_{-\infty}^{\infty}e^{-2x^2}H_0^2(x)H_m(x)H_n(x)dx,\quad m,\,n=0,\,1...,\endaligned\tag 2.5$$
The integral of the more general product 
$$\aligned &\int_{-\infty}^{\infty}e^{-2x^2}H_p(x)H_q(x)H_m(x)H_n(x)dx\\
=&(H_p(x)H_q(x)e^{-x^2},\,H_m(x)H_n(x)e^{-x^2}),\quad p,\,q,\,m,\,n=0,\,1...\endaligned\tag 2.6$$
is done in \cite{W}.

We use generating functions of Hermite polynomials to find $a_\ell$ of (1.23) as follows.
Since 
$$\align&e^{2tx-t^2}=\sum_{n=0}^\infty\frac{t^n}{n!}H_n(x),\tag 2.7\\
&e^{2sx-s^2}=\sum_{m=0}^\infty\frac{s^m}{m!}H_m(x),\tag 2.8\endalign$$
which can be found in any mathematics handbook (cf. \cite {CFKS, T} for connections with the
Mehler formula), multiplying (2.7, 2.8), we obtain 
$$\aligned &e^{2(t+s)x-(t^2+s^2)}=\sum_{n,m}\frac{t^ns^m}{n!m!}H_n(x)H_m(x)\\
=&e^{2(\frac{t+s}{\sqrt 2})\sqrt 2 x-(\frac{t+s}{\sqrt 2})^2}\cdot e^{-\frac{1}{2}(t-s)^2}\\
=&\sum_{\ell=0}^\infty H_\ell(\sqrt 2 x)\cdot\frac{(\frac{t+s}{\sqrt 2})^\ell}{\ell!}\cdot\sum_{p=0}^{\infty}
(-1)^p\frac{(t-s)^{2p}}{2^pp!}.\endaligned\tag 2.9$$

From (1.23), we are only interested in the coefficient in front of $H_0(\sqrt 2 x)$. So we set $\ell=0$.
To obtain $a_0$, we equate the coefficient in front of $t^ns^m$. Comparing the LHS with RHS of (2.9),
$n$, $m$ must have the same parity, otherwise it is $0$. We deduce
$$\aligned a_0&=\frac{(-1)^{\frac{n-m}{2}}}{2^{\frac{m+n}{2}}(\frac{m+n}{2})!}\cdot (n+m)!,\qquad n,\,m\text{ same parity}\\
&=0\qquad\qquad\qquad\qquad\qquad\qquad\text{otherwise},\endaligned\tag 2.10$$
by taking $2p=n+m$, which is the only contributing term.

Taking into account the normalization factors in the third equation of (2.3), we then obtain
\proclaim{Lemma 2.1}
$$\aligned W_{mn}{\overset\text{def }\to =}&\int_{-\infty}^{\infty}h_0^2(x)h_m(x)h_n(x)dx\\
=&\frac{(-1)^{\frac{n-m}{2}}}{2^{m+n}\sqrt{m!n!}}\cdot\frac{(m+n)!}{(\frac{m+n}{2})!}\sqrt{\frac{\pi}{2}}\qquad \text{m, n same parity}\\
=&0\qquad\qquad\qquad\qquad\qquad\qquad\quad\text{otherwise}.\endaligned\tag 2.11$$
Let $$N=\frac{n+m}{2},\, k=\frac{n-m}{2},\tag 2.12$$
assuming $n\geq m$, without loss. When $N\gg 1$,
$$\align W_{mn}&=\left[1+\Cal O\left(\frac{1}{N}\right)\right]\frac{(-1)^{\frac{n-m}{2}}N!}{\sqrt{2 N(N+k)!(N-k)!}},\tag 2.13\\
|W_{mn}|&\leq\frac{1}{\sqrt N}e^{-k^2/{2N}}.\tag 2.14\endalign$$
\endproclaim
\demo {Proof} We only need to obtain the asymptotics in (2.13, 2.14). This is an
exercise in Stirling's formula:
$$n!=\big(\frac{n}{e}\big)^n\sqrt{2\pi n}(1+\frac{1}{12n}+\frac{1}{288n^2}+...)\tag 2.15$$
or its log version
$$\log n!=(n+\frac{1}{2})\log n-n+\log\sqrt{2\pi}+...\tag 2.16$$
Here it is more convenient to use the latter. Using (2.12, 2.16),
$$\aligned &\log\frac{(m+n)!}{2^{m+n}(\frac{m+n}{2})!}=\log\frac{(2N)!}{N!}-\log 2^{2N}\\
= &N\log N-N+\frac{1}{2}\log 2+\Cal O(N^{-1}).\endaligned\tag 2.17$$
So $$\frac{(m+n)!}{2^{m+n}(\frac{m+n}{2})!}=\left[1+\Cal O\left(\frac{1}{N}\right)\right]\frac{N!}{\sqrt {\pi N}},\tag 2.18$$ 
using (2.17). Hence 
$$W_{mn}=\left[1+\Cal O\left(\frac{1}{N}\right)\right]\frac{(-1)^{\frac{n-m}{2}}N!}{\sqrt{2 N(N+k)!(N-k)!}},\quad N\gg 1$$
which is (2.13). Using the fact that
$$n!=\sqrt{2\pi n}\left(\frac{n}{e}\right)^n e^{\lambda_n}$$
with $$\frac{1}{12n+1}<\lambda_n<\frac{1}{12n},\text { for all }n\geq 1,$$
and applying the inequalities (with $x=k/N$):
$$\phi(x){\overset\text{def }\to =}(1+x)\log (1+x)+(1-x)\log (1-x)\geq x^2$$
for all $x\in[0,1)$ and $\phi(x)\geq ax^2$ with $a>1$ for $x\in[7/10, 1)$,
we obtain (2.14). (When $x=k/N=1$, (2.14) follows by a direct computation using Stirling's formula.)
\hfill $\square$
\enddemo
From (2.14), the matrix element $W_{mn}$ has subexponential decay
$$|W_{mn}|\leq\frac{1}{\sqrt N}e^{-\frac{|m-n|^{2}}{2N}}\quad (N=\frac{m+n}{2}\gg 1)\tag 2.19$$
when $|m-n|>\sqrt N$. When $|m-n|\leq \sqrt N$,
we only have the estimate 
$$|W_{mn}|\leq\frac{1}{\sqrt N}.\tag 2.20$$
Hence $W$ is a matrix with a slowly enlarging region of size $\Cal O(\sqrt N)$ around the principal diagonal $W_{nn}$, the
$\ell^2$ norm of $W$ is of $\Cal O(1)$, but the local $\ell^2$ to $\ell^\infty$ norm is of order $\Cal O(1/\sqrt N)$. These
new features will need to be taken into account when we do the analysis in sections 3, 4.

{\noindent {\it (2.1) in the Hermite-Fourier basis}

In the Hermite-Fourier basis, $e^{-in\cdot\phi}h_j(x)$, $n\in\Bbb Z^\nu$, $\phi\in\Bbb T^\nu$, $j\in\{0,\,1...\}$, 
the Floquet Hamiltonian, which is unitarily equivalent to the $K$ defined in (1.15) is then
$$H=\text {diag }(n\cdot\omega+j+\frac{1}{2})+\frac{\delta}{2}W\otimes\Delta\tag 2.21$$
on $\ell^2(\Bbb Z^\nu\times\{0,\,1...\})$,
where $W$ is the matrix operator defined in (2.11), acting on the $j$ indices,
$j=0,\,1,\,2...$,
$\Delta$ acts on the $n$ indices, $n\in\Bbb Z^\nu$,
$$\aligned \Delta_{nn'}&=1,\quad |n-n'|_{\ell^1}=1,\\
&=0,\quad \text{otherwise}.\endaligned\tag 2.22$$

Let $\tilde H=H-1/2$ and rename $\tilde H$, $H$; let $\tilde\delta=\delta/2$ and rename $\tilde\delta$, $\delta$. We then have 
$$H{\overset\text{def }\to =}\text{diag}(n\cdot\omega+j)+\delta W\otimes\Delta\tag 2.23$$
on $\ell^2(\Bbb Z^\nu\times\{0,\,1...\})$,
with $W$, $\Delta$ defined in (2.11, 2.22). 
\bigskip
\head{\bf 3. Exponential decay of Green's functions at fixed $E$: estimates in $\theta$}\endhead
Let $H$ be the operator defined in (2.23), i.e., 
$$H=\text {diag }(n\cdot\omega+j)+\delta W\otimes\Delta\tag 3.1$$
on $\ell^2(\Bbb Z^\nu\times\{0,\,1,...\})$,
where $$n\in\Bbb Z^\nu,\,j\in \{0,\,1,...\},\,0<\delta\ll 1.\tag 3.2$$
$W$ acts on the $j$ indices, $j=0,\,1,\,2...$,
$$\aligned 
W_{jj'}=&\frac{(-1)^{\frac{j-j'}{2}}}{2^{j+j'}\sqrt{j!j'!}}\cdot\frac{(j+j')!}{(\frac{j+j'}{2})!},\qquad j,\, j' \text { same
parity,}\\ =&0\qquad\qquad\qquad\qquad\qquad\qquad\text{otherwise}.\endaligned\tag 3.3$$
Write  $$J=\frac{j+j'}{2},\, k=\frac{j-j'}{2},\,(\text {assume }j\geq j').$$
When $J\gg 1$,
$$\align |W_{jj'}|&\leq\frac{J!}{\sqrt{J(J+k)!(J-k)!}}\\
&\leq \frac{1}{\sqrt J}e^{-\frac{k^2}{2J}}\tag 3.4\endalign$$
from (2.14).
$\Delta$ acts on the $n$ indices, $n\in\Bbb Z^\nu$,
$$\aligned \Delta_{nn'}&=1,\quad |n-n'|_{\ell^1}=1,\\
&=0,\quad \text{otherwise}.\endaligned\tag 3.5$$

In view of the Theorem, our aim is to prove that on a good set of $\omega$, 
$H$ has pure point spectrum. To achieve that 
goal, we add a parameter $\theta$ ($\theta\in\Bbb R$) to $H$:
$$\aligned H(\theta)&=H+\theta\\
&=\text {diag }(n\cdot\omega+\theta +j)+\delta W\otimes\Delta.\endaligned\tag 3.6$$
We consider a sequence of finite volume Green's functions at fixed $E$:
$$G_\Lambda(\theta, E)=(H_\Lambda(\theta)-E)^{-1}, \tag 3.7$$
where $\Lambda$ are cubes in $\Bbb Z^\nu\times\{0,\,1...\}$, $\Lambda\nearrow {\Bbb Z^\nu}\times\{0,\,1...\}$
in an appropriate way, $H_\Lambda(\theta)$ are $H(\theta)$ restricted to $\Lambda$.

In this section, $E\in\Bbb R$, $\omega\in [0,2\pi)^\nu$ are fixed and we assume that $\omega$ is a Diophantine 
frequency. We do estimates in $\theta$. Specifically,
we prove (inductively) that for any $\Lambda$ large enough, away from a set of $\theta$ of small measure in
$\Bbb R$, 
$\Vert G_\Lambda(\theta, E)\Vert$ is bounded and $|G_\Lambda(\theta,E)(x,\,y)|$ has subexponential decay
for $|x-y|\sim$ linear scale of $\Lambda$. (For precise statement, see Proposition 3.10.)
As mentioned in section 1, the proof is a combination of  a non perturbative part in the $n$ direction and
a perturbatuve part in the $j$ direction. We note that this
set of bad $\theta$ depends on $E$, $\omega$, the estimate on the measure is, howevere, uniform in 
$E$ and $\omega$ for Diophantine $\omega$. 

In the next section (section 4), we eliminate the $E$ dependence by excluding double resonances 
and converting the estimate in $\theta$ into estimates in $\omega$ in the process. The conversion is possible because $\omega$, $\theta$ appear in (3.6) in the form $n\cdot\omega+\theta$.

Below we start the induction process.
To simplify notations, we extend $H$ to a linear operator on $\ell^2(\Bbb Z^{\nu+1})$ with $W_{jj'}$ as in
(3.3) for $j,\,j'\in\{0,\,1...\}$ and $W_{jj'}=0$ otherwise. 

\subheading{3.1 The initial estimate ($0^{\text {th}}$ step)}

For any subset $\Lambda\subset\Bbb Z^{\nu+1}$, we define
$$\aligned H_{\Lambda}(\theta)(n,j;n',j')&=H(\theta)(n,j;n',j'),\quad
(n,j;n',j')\in\Lambda\times\Lambda,\\
&=0,\qquad\qquad\qquad\qquad\text{otherwise}.\endaligned\tag 3.8$$
Let $\Lambda_0=[-J,J]^{\nu+1}$ for some $J>0$ to be determined. 
For a fixed $E$, we study the Green's function $G_{\Lambda_0}(\theta, E)=(H_{\Lambda_0}(\theta)-E)^{-1}$
by doing perturbation theory in $\delta$. We have
\proclaim {Lemma 3.1}
Assume $0<\delta\ll 1$. For any fixed $\sigma$, $0<\sigma<1/4$, 
there exists $J\in\Bbb N$ such that the following 
statement is satisfied. Let $\Lambda_0=[-J,J]^{\nu+1}$. There 
exists a set $\Cal B(\Lambda_0, E)$ in $\Bbb R$, with 
$$\text{mes }\Cal B(\Lambda_0, E)\leq e^{-J^{\sigma/2}},\tag 3.9$$
such that if $\theta\in\Bbb R\backslash \Cal B(\Lambda_0,E )$, then 
$$\align&\Vert G_{\Lambda_0}(\theta, E)\Vert\leq e^{J^\sigma}\tag 3.10\\
&|G_{\Lambda_0}(\theta, E)(n,j;n'j')|\leq e^{-|(n,j)-(n',j')|^{1/4}},\quad \,|(n,j)-(n',j')|>J/10.\tag
3.11\endalign$$
\endproclaim
\demo{Proof} Since this is the initial estimate, we do perturbation theory in $\delta$. Let 
$\kappa>0$. Let $\Cal B(\Lambda_0, E)$ be the set such that if $\theta\in \Cal B(\Lambda_0, E)$, then
$$|n\cdot\omega+j+\theta-E|\leq 2\kappa\tag 3.12$$
for some $(n,j)\in\Lambda_0$. Clearly
$$\aligned\text {mes }\Cal B(\Lambda_0, E)&\leq 4 |\Lambda_0|\kappa\\
&=4(2J+1)^{\nu+1}\kappa.\endaligned\tag 3.13$$
Let $$D_{\Lambda_0}(\theta){\overset\text{def }\to =}\text{diag} (n\cdot\omega+j+\theta),\, (n,j)\in\Lambda_0\tag 3.14$$
be the unperturbed diagonal operator. 
Since $\Vert \delta W\otimes\Delta\Vert_{\ell^2}\leq\Cal O(\delta)$, if $\theta\notin\Cal B(\Lambda_0, E)$,
then
$$\aligned &\Vert G_{\Lambda_0}(\theta, E)\Vert=\Vert(D_{\Lambda_0}(\theta)-E+\delta W\otimes\Delta)^{-1}\Vert\\
&\leq\kappa^{-1}\endaligned\tag 3.15$$
for $\kappa\gg\delta$. 

From the resolvent equation
$$G_{\Lambda_0}(\theta,E)(n,j;n'j')=(n\cdot\omega+j+\theta-E)^{-1}[(\delta W\otimes\Delta)
G_{\Lambda_0}(\theta, E)](n,j;n'j')\tag 3.16$$
for $(n,j)\neq (n',j')$.  Hence
$$|G_{\Lambda_0}(\theta, E)(n,j;n'j')|\leq\Cal O(1)\delta\kappa^{-2}, \text{ if } |(n,j)-(n',j')|>J/10\tag 3.17$$
Let $$\align J&=|\log\delta|\text { (hence }\delta=e^{-J}), \tag 3.18\\
\kappa&=e^{-J^{\sigma}}.\tag 3.19\endalign$$
(3.13-3.17) then imply (3.9-3.11).\hfill$\square$
\enddemo

\subheading {3.2 A Wegner estimate in $\theta$ for all scales}

We now state an apriori estimate in $\theta$ for $\Vert (H_\Lambda(\theta)-E)^{-1}\Vert$ valid for  all finite subsets $\Lambda\subset\Bbb Z^{\nu+1}$ and all $\delta$. 
This estimate will be useful in the induction process. Following the Anderson localization tradition, we call it a  Wegner estimate.
\proclaim{Proposition 3.2}
For any 
$E\in\Bbb R$, and any finite subset $\Lambda$ in $\Bbb Z^{\nu+1}$, the following estimate is satisfied
for all $\kappa>0$:
$$\text {mes }\{\theta|\text { dist }(E, \sigma(H_\Lambda(\theta)))\leq\kappa\}\leq 2|\Lambda|\kappa.\tag 3.20$$
\endproclaim
\demo {Proof} Let $\lambda_k$, $k=1,..., |\Lambda|$ be eigenvalues of $H_\Lambda(\theta=0)$. Then 
$$ \{\theta|\text { dist }(E, \sigma(H_\Lambda(\theta)))\leq\kappa\}=\cup_{k=1}^{|\Lambda|}\{\theta||E- \theta-\lambda_k|\leq\kappa\}.\tag 3.21$$
It follows that the measure of the left side is bounded by $2|\Lambda|\kappa$.  \hfill $\square$
\enddemo
Our goal now is to obtain inductively the equivalent of estimates (3.9-3.11) for larger subsets
$\Lambda$, $\Bbb Z^{\nu+1}\supset\Lambda\supset\Lambda_0$. We note that in proving Lemma 3.1,
we did perturbation theory about the diagonal operator 
$$D_{\Lambda_0}(\theta)=\text{diag }(n\cdot\omega+j+\theta)|_{(n,j)\in\Lambda_0}\tag 3.22$$ 
using the smallness of $\delta$. This was sufficient for {\it one} initial scale. For 
subsequent scales, however, we need more detailed information on the spectrum of $H_{\Lambda}$.

\subheading{3.3 Local spectral property of $H_\Lambda$}

We assume that $\omega$ is Diophantine, i.e., 
$\exists\, c>0,\,\alpha>2\nu$, such that 
$$|n\cdot\omega+j|\geq\frac{c}{|n|^\alpha}\tag 3.23$$
for all $n\in\Bbb Z^\nu\backslash\{0\}$, all $j\in\Bbb Z$. In this subsection,
we make statements which hold for {\it any} fixed $\theta$. We look at finite subsets $\Lambda\subset \Bbb
Z^{\nu+1}$, such that $j\neq 0$, if $(n,j)\in\Lambda$. 
When $j\gg 1$, this is the perturbative region. 
\proclaim{Proposition 3.3} Assume $\omega$ satisfies (3.23). Let $0<\beta'<1/5\alpha$ and $3/4<\beta<1$.
Let $\Lambda$ be a rectangle centered at ($N,\,2L)\in\Bbb Z^\nu\times\Bbb Z$, where $\Bbb Z^\nu$
is identified with $\Bbb Z^\nu\times\{0\}$:
$$\Lambda=(N,\,2L)+[-L^{\beta'},L^{\beta'}]^{\nu}\times [-L^{\beta},L^{\beta}] \subset\Bbb Z^{\nu+1}.\tag 3.24$$ 
Assume $$|L|\gg 1.\tag 3.25$$
For any fixed $\theta\in\Bbb R$, the eigenvalues $\lambda_{n,j}(\theta)$ of $H_\Lambda(\theta)$ satisfy 
$$|\lambda_{n,j}(\theta)-\lambda_{n',j'}(\theta)|>\frac{1}{L^{1/5}}\quad (n,\,j)\neq (n'\,j');\tag
3.26$$ the eigenfunctions $\phi_{n,j}$ may be chosen such that 
$$\Vert \phi_{n,j}-\delta_{n,j}\Vert<\frac{1}{L^{1/20}}.\tag 3.27$$
\endproclaim

\noindent{\it Remark.} It is
crucial to note that the estimates (3.26, 3.27) are {\it independent} of $\theta$ and the specific $\omega$
satisfying (3.23). In Lemma 3.5, we exploit further
the consequences of (3.26).
\demo{Proof of Proposition 3.3}
Let $(n,\,j)$, $(n',\,j')\in\Lambda$, $(n,\,j)\neq (n',\,j')$. Then {\it for all} $\theta$, the difference
of the diagonal elements
$$|n\cdot\omega+j+\theta-(n'\cdot\omega+j'+\theta)|=|(n-n')\cdot\omega+(j-j')|\geq\frac{c}{L^{\beta'\alpha}}>
\frac{1}{L^{1/5}}\text{ for }0<\beta'<1/5\alpha,\, L\gg 1,\tag 3.28$$
from (3.23). 
Use as approximate eigenfunctions $\delta_{n,j}$ with approximate eigenvalues
$\tilde \lambda_{n,j}(\theta)=n\cdot\omega+j+\theta$,
$(n,j)\in\Lambda$, and let
$$(H_{\Lambda}-\tilde \lambda_{n,j})\delta_{n,j}=\psi.\tag 3.29$$
Then from (3.7),
$$\aligned \psi(n',j')&=\delta W_{j'j},\quad |n-n'|=1,\,(n',\,j')\in\Lambda\\
&=0,\quad\qquad \text{otherwise},\endaligned\tag 3.30$$
and $$\Vert\psi\Vert_{\ell^2}=\Cal O(|L|^{{-1/4}})\tag 3.31$$
from (3.4). 

Equations (3.28-3.31) imply that $\tilde \lambda_{n,j}(\theta)=n\cdot\omega+j+\theta$ is an approximate eigenvalue of
$H_\Lambda$ to $\Cal O (|L|^{{-1/4}})$. This can be seen as follows. Let $\lambda=\tilde\lambda_{n,j}(\theta)$.
Assume
$$\text{ dist }(\lambda,\sigma(H_{\Lambda}(\theta)))>\Cal O(|L|^{{-1/4}}).$$
Take any $f$, $\Vert f\Vert_{\ell^2(\Lambda)}=1$, $f=\sum c_m\psi_m$, $\sum|c_m|^2=1$, where
$\psi_m$ are eigenfunctions of $H_{\Lambda}$. Then 
$$(H_{\Lambda}-\lambda)f=\sum(E_m-\lambda)c_m\psi_m,$$
where $E_m$ is the corresponding eigenvalue for $\psi_m$. 
Hence
$$\Vert(H_{\Lambda}-\lambda)f\Vert=\sqrt {\sum (E_m-\lambda)^2|c_m|^2}>\Cal O(|L|^{{-1/4}}).\tag 3.32$$
This is a contradiction if $f=\delta_{n,j}$.

Hence (3.28) gives (3.26). (3.27) follows from (3.26, 3.31) and standard perturbation theory, see e.g., \cite{Ka}.
\hfill $\square$
\enddemo
\smallskip
\subheading {3.4. The first iteration ($1^{\text {st}}$ step)}

We now increase the scale from $J$ to $J^C$, where $C>1$ (independent of $\delta$) is the geometric expansion factor, which will be specified in section 3.5. Recall that
$J$ is the initial scale, large enough so that (3.18) holds. Hence on $\Lambda_0=[-J,J]^{\nu+1}$, (3.9-3.11) hold.

Let $$\Lambda\overset\text{def }\to =[-J^C,J^C]^{\nu+1}.\tag 3.33$$ 
Our aim is to prove the analogue of (3.9-3.11) when $\Lambda_0$ is replaced by $\Lambda$. The
general strategy in going from scale $J$ to scale $J^C$ is to distinguish the region near
$j=0$, where we use the estimate from the previous scale, here (3.9-3.11) and non perturbative 
arguments and the region 
away from $j=0$, where we use Proposition 3.3 and perturbation theory. The general iteration strategy
here is similar to that in \cite{BW1}.

Toward that end, we define
$$\Cal T\overset\text{def }\to=\{(n,j)\in\Lambda|\,|j|\leq 2J-1\}.\tag 3.34$$
Let $$\Lambda_0(n,j)=(n,j)+[-J,J]^{\nu+1}, \quad (n,j)\in\Cal T\tag 3.35$$
be cubes of the previous scale. Let 
$$\Lambda_*(n,2j)=(n,2j)+[-j^{\beta'},j^{\beta'}]^{\nu}\times[-j^\beta,j^\beta]\quad
(0<\beta'<1/5\alpha,\,3/4<\beta<1),
\quad (n,2j)\in\Lambda\backslash\Cal T,\tag 3.36$$
be cubes of type (3.24). We cover $\Lambda$ with $\Lambda_0$, $\Lambda_*$ cubes, i.e., $\Cal T$ with
$\Lambda_0$, $\Lambda\backslash\Cal T$ with $\Lambda_*$.
$G_\Lambda$ are then obtained by using resolvent equation, (3.9-3.11), Propositions 3.2 and 3.3. We implement
this strategy in our first iteration. This iteration is special, as in the tube region $\Cal T$,
we use smallness of $\delta$, cf. (3.18).

We need the following notion of pairwise disjointness. Let $S_k$, 
$k=1,...,K$ be finite sets , $S_k\neq S_{k'}$ if $k\neq k'$. Let $S=\{S_k\}_{k=1}^K$.  If $S_k\cap S_{k'}\neq\emptyset$, $\forall k\neq k'$, then we say  there is {\it 1 pairwise disjoint set } in $S$. More 
generally, if $\exists I_1,\,I_2,...,I_P$, $I_p\cap I_{p'}=\emptyset$, if $p\neq p'$, $\{I_p\}_{p=1}^P=\{1,\,2,...,K\}$ such that $S_k\cap S_{k'}\neq\emptyset$ if and only if $k$, $k'\in I_p$ for some $p$.  Then
we say there are $P$ pairwise disjoint sets in $S$.

\proclaim{Lemma 3.4} Let $\Gamma_0$ be a covering of $\Cal T$ with $\Lambda_0$
cubes defined in (3.35). Assume $\omega$ is Diophantine satisfying (3.23). Fix $E$, $\sigma\,(0<\sigma<1/4)$ as in Lemma 3.1. For 
all $\theta$, there exists at most {\bf 1} pairwise disjoint $\tilde \Lambda\in\Gamma_0$,
such that 
$$\text{ dist }(E,\,\sigma(H_{\tilde\Lambda}(\theta)))\leq e^{-J^{\sigma}}.\tag 3.37$$
Moreover if $\Lambda_0\in\Gamma_0$ and $\Lambda_0\cap\tilde\Lambda=\emptyset$,
$$\text{ dist }(E,\,\sigma(H_{\Lambda_0}(\theta)))\geq \frac{c}{J^{C\alpha}}.\tag 3.38$$
\endproclaim
\demo{Proof}
Let $\Lambda_0\in\Gamma_0$. If $|\lambda_n-E|<e^{-J^\sigma}$ for a $\lambda_n\in\sigma(H_{\Lambda_0}(\theta)$,
then since $\Vert\delta W\otimes\Delta\Vert\leq C\delta=Ce^{-J}$,  
$$|n\cdot\omega+j+\theta-E|\leq e^{-J^{\sigma}}+Ce^{-J}\leq 2 e^{-J^{\sigma}}$$ 
for some $(n,j)\in\Lambda_0$. Thus if both $\Lambda_0$ and $\tilde\Lambda$ satisfy (3.37), then 
$$|(n-n')\cdot\omega+j-j'|\leq |n\cdot\omega+j+\theta-E|+|n'\cdot\omega+j'+\theta-E|\leq 4 e^{-J^{\sigma}}\tag 3.39$$
for some $(n,j)\in\Lambda_0$ and $(n',j')\in\tilde\Lambda$. This implies $\Lambda_0\cap\tilde\Lambda\neq\emptyset$
for large $J$, since otherwise the left side is larger than $c/(2\nu J^{C})^\alpha$
by the Diophantine condition (3.23), which is a contradiction if $C\ll J^\sigma/\log J$. 

If $\Lambda_0\in\Gamma_0$ is disjoint from $\tilde\Lambda\in\Gamma_0$,  then for any $(n,j)\in\Lambda_0$ and
for some $(n', j')\in\tilde\Lambda$, 
$$ |n\cdot\omega+j+\theta-E|\geq |(n-n')\cdot\omega+j-j'|-|n'\cdot\omega+j'+\theta-E|\geq \frac{c}{(2\nu J^{C})^\alpha}-2e^{-J^{\sigma}}.\tag 3.40$$
Since $\Vert\delta W\otimes\Delta\Vert\leq C\delta=Ce^{-J}$, (3.38) follows.
\hfill$\square$
\enddemo
If $\Lambda_0\cap \tilde\Lambda=\emptyset$, then (3.10, 3.11) are available. Let $\Gamma_*$ be
a covering of $\Lambda\backslash\Cal T$. For $\Lambda_*\in\Gamma_*$, we need 
\proclaim{Lemma 3.5} Let $0<\beta'<1/5\alpha$ and $3/4<\beta<1$. Let $\Lambda_*$ be the rectangle: 
$$\Lambda_*(n,2j)=(n,2j)+[-j^{\beta'},j^{\beta'}]^{\nu}\times[-j^\beta,j^\beta],
\quad (n,2j)\in\Lambda\backslash\Cal T.\tag 3.41$$
For a fixed $E\in\Bbb R$,  there exists $W\subset\Bbb R$, with $\text {mes } W\leq C|\Lambda_*|e^{-J^{\beta'/2}}$ 
such that for $\theta\in\Bbb R\backslash W$, $J\gg 1$, 
$$\align &|[H_{\Lambda_*}(\theta)-E]^{-1}(n',j';n'',j'')|\leq
e^{-\frac{1}{2}[|n'-n''|+\frac{|j'-j''|}{j^{1/2}}]},\\ &\qquad\qquad\qquad\qquad (|n'-n''|>j^{\beta'}/10,\text{ or }|j'-j''|\geq
j^\beta/10)\tag 3.42\\
&\Vert [H_{\Lambda_*}(\theta)-E]^{-1}\Vert\leq e^{J^{\beta'/2}}.\tag 3.43\endalign$$
\endproclaim
\demo{Proof}
For any fixed $\theta$, $E$, there exists at most $1$ bad site $b=(n',j')\in\Lambda_*$, such that 
$$|n'\cdot\omega+j'+\theta-E|\leq \frac{1}{j^{1/5}}.$$
This is because if there were $(n',j')$, $(n'',j'')\in\Lambda_*$,  $(n',j')\neq (n'',j'')$ such that 
$$|n'\cdot\omega+j'+\theta-E|\leq \frac{1}{j^{1/5}},$$
$$|n''\cdot\omega+j''+\theta-E|\leq \frac{1}{j^{1/5}},$$
then 
$$|(n'-n'')\cdot\omega+(j'-j'')|\leq\frac{2}{j^{1/5}},$$
which contradicts the Diophantine condition (3.23) on $\omega$:
$$\Vert (n'-n'')\cdot\omega\Vert_{\Bbb T}\geq\frac{c}{|n'-n''|^\alpha}\geq\frac{c}{j^{\beta'\alpha}}\gg\frac{2}{j^{1/5}}\tag 3.44$$
for $j\geq J\gg1$, using (3.41).

To obtain (3.42, 3.43), we first consider $H_{\Lambda_*\backslash b}(\theta)$ and make estimates on 
$[H_{\Lambda_*\backslash b}(\theta)-E]^{-1}(n',j';n'',j'')$. We then use Proposition 3.2 and the resolvent
equation to obtain (3.42). We prove below that 
$$\Vert [H_{\Lambda_*\backslash b}(\theta)-E]^{-1}\Vert_{\ell^2\to\ell^2}<2j^{1/5},\tag 3.45$$
for $j\gg 1$.
To obtain estimates on the matrix elements, we use weighted $\ell^2$ space and show that 
$[H_{\Lambda_*\backslash b}(\theta)-E]^{-1}$ remain bounded. (3.45) is then the special case when the weight equals
to $1$.

Toward that end, for any $p\in\Bbb Z^{\nu+1}$, we define $|p|_\nu=|(p_1,...,p_\nu,0)|$ and $|p|_1=|(0,...,0,p_{\nu+1})|$.
Let $\rho=\{\rho_a\}_{a\in\Bbb Z^{\nu+1}}$ be a family of
weights, such that 
$$\rho_a(p)=e^{|p-a|_\nu+\frac{|p-a|_1}{j^{1/2}}}\quad (a\in\Bbb Z^{\nu+1}),\quad \forall\, p\in\Bbb
Z^{\nu+1},\tag 3.46$$
where $j$ is as in (3.41).
Hence $\forall\,p,\,q\in\Bbb Z^{\nu+1}$,
$$\aligned \rho_a^{-1}(p)\rho_a(q)&=e^{-|p-a|_\nu-\frac{|p-a|_1}{j^{1/2}}}\cdot
e^{|q-a|_\nu+\frac{|q-a|_1}{j^{1/2}}}\\ &\leq e^{|p-q|_\nu+\frac{|p-q|_1}{j^{1/2}}},\endaligned\tag
3.47$$ for all $\rho_a\in\rho$.

Let $D$ be the diagonal part in (3.8). To arrive at (3.42), we consider the deformed operator
$\tilde H_{\Lambda_*\backslash b}(\theta)$:
$$\aligned \tilde H_{\Lambda_*\backslash b}(\theta)(p,q)\overset\text{def }\to =&[\rho_a^{-1}
H_{\Lambda_*\backslash b}(\theta)\rho_a](p,q)\\
=&D_{pq}+\delta\rho_a^{-1}(p)\rho_a(q)(W\otimes\Delta)_{pq},\\
=&D_{pq}+\delta(\tilde W\otimes\tilde\Delta)_{pq}\endaligned\tag 3.48$$
where $p,\,q\in\Lambda_*\backslash b$ and 
$$(\tilde W\otimes\tilde\Delta)_{pq}{\overset\text{def }\to =}\rho_a^{-1}(p)\rho_a(q)(W\otimes\Delta)_{pq}.\tag 3.49$$

To prove boundedness of $[\tilde H_{\Lambda_*\backslash b}(\theta)-E]^{-1}$, we use resolvent
series and perturb about the diagonal. We have formally
$$\aligned [\tilde H_{\Lambda_*\backslash b}(\theta)-E]^{-1}=&(D-E)^{-1}+\delta(D-E)^{-1}(\tilde W\otimes\tilde\Delta)(D-E)^{-1}\\
+&\delta^2(D-E)^{-1}(\tilde W\otimes\tilde\Delta)(D-E)^{-1}(\tilde W\otimes\tilde\Delta)(D-E)^{-1}\\
+&...\endaligned\tag 3.50$$
Let $$\Cal W{\overset\text{def }\to =}(\tilde W\otimes\tilde\Delta)(D-E)^{-1}.$$ Then 
$$\aligned \Cal W(n',j';n'',j'')&=\tilde W_{j'j''}\tilde\Delta_{n'n''}(n''\cdot\omega+j''+\theta-E)^{-1},\quad |n'-n''|=1\\
&=0\quad \text{otherwise}.\endaligned\tag 3.51$$
We bound $\Cal W$ using Schur's lemma \cite{Ka}:
$$\Vert\Cal W\Vert\leq(\sup_{n',j'}\sum_{n'',j''}|\Cal W(n',j';n'',j'')|)^{1/2}(\sup_{n'',j''}\sum_{n',j'}|\Cal W(n',j';n'',j'')|)^{1/2},\tag 3.52$$
where $(n',j')$, $(n'', j'')\in\Lambda_*\backslash b$.
$$\sum_{n'',j''}|\Cal W(n',j';n'',j'')|=\sum_{|n''-n'|=1}\sum_{|j''-2j|\leq j^\beta}|\tilde W_{j'j''}\tilde\Delta_{n'n''}(n''\cdot\omega+j''+\theta-E)^{-1}|,$$
$(n', j')\in\Lambda_*\backslash b$.
From (3.47, 3.4)
$$\aligned |\tilde W_{j'j''}\tilde\Delta_{n'n''}|&\leq \frac{e}{\sqrt j} e^{[-\frac{(j'-j'')^2}{2(j'+j'')}+\frac{|j'-j''|}{j^{1/2}}]}\\
&\leq \frac{e^{200}}{\sqrt j} e^{-\frac{(j'-j'')^2}{10j}},\endaligned\tag 3.53$$
where we used $j'+j''\leq 4j+2j^\beta\leq 9j/2$ for $j'$, $j''\in\Lambda_*\backslash b$.
So $$\aligned \sum_{n'',j''}|\Cal W(n',j';n'',j'')|&\leq  \frac{e^{200}}{\sqrt j}\sum_{|n''-n'|=1}\sum_{|j''-2j|\leq j^\beta}
|n''\cdot\omega+j''+\theta-E|^{-1}\\
&\leq \frac{\Cal O(1)}{\sqrt j}(j^{1/5}+\log j).\endaligned\tag 3.54$$
Let $a_{j''}{\overset\text{def }\to =}n''\cdot\omega+j''+\theta-E.$ To arrive at (3.54), we used the fact that 
$|a_{j''}|=|n''\cdot\omega+j''+\theta-E|> j^{-1/5}$ for all $(n'',j'')\in\Lambda_*\backslash b$
and that $a_p-a_q=p-q$, for all $p$, $q$. 

Using (3.53), we have 
$$\sum_{n',j'}|\Cal W(n',j';n'',j'')|\leq\Cal O(1)j^{1/5}.\tag 3.55$$
Substituting (3.54, 3.55) into (3.52), we have
$$\aligned\Vert\Cal W\Vert&\leq\Cal O(1)\big(\frac{j^{2/5}}{j^{1/2}}\big)^{1/2}\\
&< \Cal O(1)j^{-1/20}\quad (j\gg1).\endaligned$$

So the Neumann series in (3.50) is norm convergent:
$$\aligned \Vert [\tilde H_{\Lambda_*\backslash b}(\theta)-E]^{-1}\Vert&\leq\Vert(D-E)^{-1}\Vert(1+\delta\Vert \Cal W\Vert+\delta^2\Vert \Cal W\Vert^2+...)\\
&<2\Vert(D-E)^{-1}\Vert\\
&<2j^{1/5},\endaligned$$
which is equivalent to 
$$\Vert \rho_a(H_{\Lambda_*\backslash b}(\theta)-E)^{-1}\rho_a^{-1}\Vert< 2j^{1/5}.\tag 3.56$$
For each pair $p,\,q\in\Lambda_*\backslash b$, we can then always choose $a$ so that
$$|[H_{\Lambda_*\backslash b}(\theta)-E]^{-1}(n',j';n'',j'')|\leq Cj^{1/5}
e^{-[|n'-n''|+\frac{|j'-j''|}{j^{1/2}}]}.\tag 3.57$$
Write $G_0$ for $[H_{\Lambda_*\backslash b}(\theta)-E]^{-1}$, $G$ for $[H_{\Lambda_*}(\theta)-E]^{-1}$.
From the resolvent equation:
$$G=G_0+G_0H_bG_0+G_0H_bGH_bG_0,$$
where $H_b=H_{\Lambda_*}-H_{\Lambda_*\backslash b}$. The matrix element 
$H_b(p,q)=0$, unless $p=b$ or $q=b$. Using this, (3.57), the Wegner estimate (3.20) on 
$G$ with $\kappa= e^{-J^{\beta'/2}}$ ($0<\beta'<1/5\alpha$) and self-adjointness, we obtain (3.42, 3.43)
\hfill$\square$
\enddemo

We now write the estimate at scale $J^C$. Assume $\delta$, $J$ satisfying (3.18), so that Lemma 3.1 holds. Let
$J_1=J^C$, $C>1$, the same geometric expansion factor as before. Let $\Lambda=[-J_1,J_1]^{\nu+1}$.
We have
\proclaim{Lemma 3.6}
Assume $\omega$ is Diophantine satisfying (3.23) and $0<\delta\ll 1$ is the same as in Lemma 3.1. 
For any fixed $\sigma$, $0<\sigma<1/5\alpha$, there exists $\Cal B(\Lambda, E)$ in $\Bbb R$, with
$$\text{mes }\Cal B(\Lambda, E)\leq e^{-J_1^{\sigma/2}},\tag 3.58$$
such that if $\theta\in\Bbb R\backslash \Cal B(\Lambda,E)$, then 
$$\align&\Vert G_{\Lambda}(\theta, E)\Vert\leq e^{J_1^{\sigma}}\tag 3.59\\
&|G_{\Lambda}(\theta, E)(n,j;n'j')|\leq e^{-|(n,j)-(n',j')|^{1/4}},\\
&\qquad\qquad \text{for all } (n,j),\,(n'j')\,\text{such that }|(n,j)-(n',j')|>J_1/10,\tag 3.60\endalign$$
provided the expansion factor $C$ satisfies $1<C<\beta'/\sigma$, $0<\beta'<1/5\alpha$ is
as in (3.41).
\endproclaim
So we have the same estimate as in Lemma 3.1 at the larger scale $J_1=J^C$ ($C>1$). 
\demo{Proof}
This is similar to the proof of Lemma 2.4 in \cite {BW1}. So we summarize the main steps.
We prove (3.59, 3.60) using the resolvent equation and cover $\Lambda$ with cubes of types $\Lambda_0$,
$\Lambda_*$ defined in (3.35, 3.36), i.e., $\Cal T$, defined in (3.34), with a covering $\Gamma_0$ of $\Lambda_0$'s
and $\Lambda\backslash \Cal T$, a covering $\Gamma_*$ of $\Lambda_*$'s.  From Lemma
3.4, for all fixed $E$, all $\theta$, there exists at most $1$ pairwise disjoint
$\Lambda_0\in\Gamma_0$ on which (3.10, 3.11) do not hold. We use Proposition 3.2 on this $\Lambda_0$, (3.10, 3.11) on
all other $\Lambda_0$.

For a given $\Lambda_*$, let $W_{\Lambda_*}$ be the set such that (3.42, 3.43) hold if 
$\theta\in\Bbb R\backslash W_{\Lambda_*}$. Let $\Cal W=\cup W_{\Lambda_*}$, where the union is over 
all possible $\Lambda_*$ with centers in $\Lambda\backslash\Cal T$.
$$\text{mes } \Cal W\leq\Cal O(1)J_1^{2(\nu+1)}e^{-J^{\beta'/2}}\leq 
e^{-J_1^{\frac{\beta'}{2C}}}\cdot J_1^{2(\nu+1)},$$
where the first term is an upper bound on the number of possible $\Lambda_*$ with centers in 
$\Lambda\backslash\Cal T$ multiplied by  the volume of $\Lambda_*$.
For $\theta\in\Bbb R\backslash\Cal W$, we use (3.42, 3.43) on $\Lambda_*$. 

Using the resolvent equation, combining (3.10, 3.11, 3.42, 3.43) and (3.20) with $\kappa= e^{-J_1^{\sigma}}$
($0<\sigma<1/5\alpha$) on the only bad $\Lambda_0$, we obtain (3.59, 3.60). (For more general iterations
using the resolvent equation, cf. proof of Lemma 3.8, in particular (3.72).) 
Combining the measure estimate from (3.20)
with $\kappa= e^{-J_1^{\sigma}}$ and the above measure estimate on $\Cal W$, we obtain (3.58),
provided $C<\beta'/\sigma$.
\hfill$\square$
\enddemo

For the induction process to follow, it is convenient to define the following. For any fixed $\sigma$, $0<\sigma<1/5\alpha$ and any given box
$\Lambda\subset\Bbb Z^{\nu+1}$ of side length $2J+1$, we say $G_\Lambda(\theta, E)$ at
fixed ($\theta$, $E$) is {\it good} if
$$\align&\Vert G_{\Lambda}(\theta, E)\Vert\leq e^{J^{\sigma}}\\
&|G_{\Lambda}(\theta, E)(n,j;n'j')|\leq e^{-|(n,j)-(n',j')|^{1/4}},\\
&\qquad\qquad \forall (n,j),\,(n'j')\,\text{such that }|(n,j)-(n',j')|>J/10.\tag 3.61 \endalign$$
Otherwise, it is {\it bad}. We also define
$$\Cal W(\Lambda){\overset\text{def }\to =}\bigcup W_{\Lambda_*},\tag 3.62$$
where the union is over all possible $\Lambda_*$ of the form (3.36) with centers in $\Lambda\backslash\Cal T$
and $\Cal T$ is as defined in (3.34), so that for 
$\theta\in\Bbb R\backslash\Cal W(\Lambda)$, (3.42, 3.43) are valid for all $\Lambda_*$ with centers in
$\Lambda\backslash\Cal T$. We say $\Lambda_*$ is good, if (3.42, 3.43) hold. Otherwise $\Lambda_*$ is bad.
\smallskip
\subheading{3.5 A large deviation estimate in $\theta$ for the Green's functions at fixed $E$ at all scales}

We now increase the scale from $J_1$ to $J_1^C$ ($1<C<\beta'/\sigma$, the geometric expansion factor will be 
determined here). Our task is again to derive estimates (3.58-3.60) for the cube $[-J_1^C,J_1^C]^{\nu+1}$
starting from the estimates (3.58-3.60) for the cube $[-J_1,J_1]^{\nu+1}$. 

For simplicity of notation, we rename $J_1$, $J$ and $[-J_1^C,J_1^C]^{\nu+1}$, $\Lambda$ in this 
section. As in the first iteration, we distinguish the
tube region $\Cal T$, defined as in (3.34) with the new $J$. We cover $\Cal T$ with 
$\Lambda_0$'s defined in (3.35) with the new $J$, and $\Lambda\backslash\Cal T$,
$\Lambda_*$ defined in (3.36). We note that $\Lambda_0$ are at scale $J$ with centers in $\Cal T$,
while $\Lambda_*$ are at scales from $J^{\beta'}$ to $J^{C\beta}$ ($0<\beta'<1/5\alpha,\,3/4<\beta< 1$)
with centers away from $\Cal T$. 

As in the first iteration, we use the resolvent equation to obtain estimates on $G_\Lambda$ from estimates
on $G_{\Lambda_0}$ and $G_{\Lambda_*}$. Let $\Gamma_0$ be a covering of $\Cal T$ and $\Gamma_*$
of $\Lambda\backslash\Cal T$. For
$\theta\in\Bbb R\backslash\Cal W(\Lambda)$, (3.42, 3.43) are valid on  {\it all} $\Lambda_*\in\Gamma_*$. So for
any  {\it fixed} $\theta\in\Bbb R\backslash\Cal W(\Lambda)$, we only need to control the number of 
pairwise disjoint bad $\Lambda_0$ boxes on which estimate (3.61) is not available. In particular, we need the
number of such bad boxes to be $\ll J^C$, the linear scale of the box $\Lambda$. (This is intuitively clear,
as otherwise without further detail on the location of the bad boxes, we could not accumulate decay at the
linear scale as in (3.60).) Recall that for the first iteration, there is at most $1$ such (pairwise
disjoint) bad box.

\proclaim{Lemma 3.7}
Assume $\omega$ is Diophantine, satisfying (3.23) and Lemma 3.6 is valid on cubes
$$\Lambda_0(0,j)=(0,j)+[-J,J]^{\nu+1},\,\forall\,j\in[-(2J-1),2J-1].\tag 3.63$$
Then for all {\it fixed} $\theta$
$$\aligned &\#\{(n,j)\in\Cal T|\Lambda_0(n,j)\text { is a bad box }\}\\
\leq&\Cal O(1)J^{5(\nu+1)}=(J^C)^{1-}\ll J^C,\endaligned\tag 3.64$$
by choosing $5(\nu+1)<C\ll J^{\sigma/2}$ ($0<\sigma<1/5\alpha$).
\endproclaim
\demo {Proof}
Write $\Lambda_0$ for $\Lambda_0(0,j)$. We first replace the estimate 
$\Vert G_{\Lambda_0}(\theta, E)\Vert_{\ell^2\to\ell^2}\leq e^{J_1^{\sigma}}$ in (3.59)
by the estimate on the Hilbert-Schmidt norm:
$$\Vert G_{\Lambda_0}(\theta, E)\Vert_{\text {HS}}\leq e^{J_1^{\sigma}}.\tag 3.59'$$
This leaves the measure estimate in (3.58) unchanged for $J_1\gg 1$ (cf. proof of Lemma 3.6).

Define
$$
\Cal A\overset {\text {def}}\to= \bigcup\limits_{j\in [-(2J-1),\, 2J-1]}
\ \Cal B(\Lambda_0 (0,j)).\tag 3.65
$$
Since the conditions on the Green's function in (3.59', 3.60) can be rewritten as polynomial inequalities in 
$\theta$ by
using Cramer's rule, $\Cal A$ is semi-algebraic of total degree less
than
$$
\align
&(2J+1)^{2(\nu+1)} \cdot (2J+1)^{2(\nu+1)}\cdot (4J+1)\\
&\leq \Cal O_{\nu}(1) J^{5(\nu+1)},\tag 3.66
\endalign
$$
where the first factor is an upper bound of the degree of polynomial for each entry of the matrix
$G_{\Lambda_0}(\theta, E)$, the second is an upperbound on the $\#$ of entries of each $G_{\Lambda_0}(\theta, E)$ plus the one for the
Hilbert-Schmidt norm, the third is the $\#$ of different matrices $G_{\Lambda_0}$'s. 
For more details, cf. the proof 
of Lemma 2.6 in \cite{BW1}.
$\Cal A$ is therefore the union of at  most $\Cal O_{\nu}(1) J^{5(\nu+1)}$ intervals in $\Bbb R$ by using
Theorem 1 in [Ba] (see also [BGS], where the special case we need is restated as Theorem 7.3).

For any fixed $\theta\in \Bbb R$, let
$$
I=\{n\in [-J^C, J^C]^{\nu}\big| n\cdot\omega+\theta\in \Cal A\}.\tag 3.67
$$
Then $\Lambda_0(n,j)$ is a bad box if and only if $n\in I$ and, therefore, $n$ is in one of the intervals of 
$\Cal A$. But each interval does not contain two such $n$ if 
$\omega$ satisfies (3.23) and
$|I|\leq \Cal O_{\nu}(1) J^{5(\nu+1)}$ by virtue of (3.58). This is because for Diophantine
$\omega$ satisfying (3.23), if there exist
$n,\,n'\in[-J^C,J^C]^\nu$, $n\neq n'$, then 
$$
|(n-n')\cdot\omega\big|\geq \frac c{(2J^C)^\alpha}\gg e^{-J^{\sigma/2}}.\tag 3.68
$$
Hence each interval can contain at most 1 integer point in $[-J^C, J^C]^{\nu}$.
We therefore obtain (3.64).\hfill $\square$
\enddemo

For any fixed $\theta\in\Bbb R\backslash\Cal W (\Lambda)$, $\Cal W(\Lambda)$ defined as in (3.62),
the only bad boxes are of type $\Lambda_0$. Lemma 3.7 shows that there are only few (of order ($(J^C)^{1-}$)
bad $\Lambda_0$ boxes in $\Lambda$. The following iteration lemma will enable us to obtain estimates 
(3.58-3.60) for $G_\Lambda$ at scale $J^C$.

\proclaim{Lemma 3.8} Fix $b\in(0, 1/8)$ and assume $\tau$ satisfies $3/4+b<\tau<1-b$. Suppose $M$, $N$ 
are integers satisfying $$N^\tau\leq M\leq 2N^\tau.\tag 3.69$$
Let $\Lambda=[-N, N]^{\nu+1}$. Assume for all $\bar \Lambda\subseteq\Lambda$ with diameter $L$, the 
Green's function $G_{\bar\Lambda}(E)=(H_{\bar\Lambda}-E)^{-1}$ at energy $E$ satisfies 
$$\Vert G_{\bar\Lambda}\Vert\leq e^{L^b}. \tag 3.69'$$

Let $\Lambda'$ be cubes of side length $2M$. We say that $\Lambda'$ is good if in addition to (3.69'), the Green's 
function exhibits off-diagonal decay:
$$|G_{\Lambda'}(E)(x, y)|\leq e^{-|x-y|^{1/4}}$$
for all $x, y\in\Lambda'$ satisfying $|x-y|>M/10$.  Otherwise $\Lambda'$ is bad. Assume for any family $\Cal F$
of pair-wise disjoint bad $\Lambda'$ cubes in $\Lambda$,
$$\# \Cal F\leq N^b.$$
Under these assumptions, one has 
$$|G_\Lambda(E)(x, y)|\leq e^{-|x-y|^{1/4}}$$
for all $x, y\in\Lambda$ satisfying $|x-y|\geq N/10$, provided $N$ is sufficiently large, i.e., $N\geq N_0(b,\tau)$.
\endproclaim 
\demo{Proof} The proof is similar to the proof of Lemma 2.4 in \cite{BGS}. As we will see from (3.73, 3.74), because of the 
conditions on $b, \tau$, it only needs a one step iteration.

To estimate $G_\Lambda(E)(x, y)$, $x, y\in\Lambda$, $|x-y|\geq N/10$, let $Q$ be cubes of side length $4M$, we make an exhaustion $\{S_i(x)\}_{i=0}^\ell$
of $\Lambda$ of width $2M$ centered at $x$ as follows:
$$\aligned &S_{-1}(x)\overset\text{def }\to=\emptyset,\\
&S_{0}(x)\overset\text{def }\to=\Lambda'(x)\cap\Lambda,\\
&S_{i}(x)\overset\text{def }\to=\cup_{y\in S_{i-1}}Q(y)\cap\Lambda,\endaligned$$
for $1\leq i\leq\ell$, where $\ell$ is maximal such that $S_\ell(x)\neq\Lambda$.

We say an annulus $A_i=S_i(x)\backslash S_{i-1}(x)$ is good if $A_i\cap\Cal F=\emptyset$. Let $A_i(x)$, 
$A_{i+1}(x),...,A_{i+s}(x)$ be adjacent good annuli and define 
$$U=\cup_{k=i}^{i+s} A_k(x).$$
Let $\partial_*S_{-1}(x)=\{x\}$ and 
$$\partial_*S_{j}(x)=\{y\in S_j(x)|\exists z\in\Lambda\backslash S_j(x), |y-z|=1\}$$
for $j\geq 0$.

By construction $$\text {dist }(\partial_*S_{i-1}, \partial_*S_{i+s})=2M(s+1)\quad (s\geq 0).$$
For any subset $B\subseteq\Lambda$, let $H_B$ be defined as in (3.8). Let 
$$\Gamma_B \overset\text{def }\to=H_\Lambda-(H_B\oplus H_{\Lambda\backslash B}).$$
Assume $x\in B$, from the resolvent equation,
$$\aligned G_\Lambda(E)(x,&y)=G_{B}(E)(x,y)\\
&+\sum_{\Sb z\in B\\ z'\in \Lambda\backslash B\endSb}G_{B}(E)(x,z)\Gamma_B
(z,z')G_\Lambda(E)(z',y).\endaligned \tag 3.70$$
The proof follows the same line of arguments as in the proof of Lemma 2.4 in \cite{BGS} by iterating (3.70). There are 
two modifications:

\item\item{$\bullet$} In view of $\Gamma_B$, with matrix elements
$\Gamma_B(z, z')=\delta W(j,j')\Delta(n,n')$, $z=(n,j)\in B$, $z'=(n',j')\in\Lambda\backslash B$,
where $\Delta(n, n')$ as defined in (2.22) and 
$$|W(j,j')|\leq e^{-\frac{(j-j')^2}{2(j+j')}}\leq e^{-\frac{(j-j')^2}{4N}}$$
from (2.14) for  $1\ll j+j'\leq 2N$, 
we define
$$\aligned &\partial S_i^{(1)}(x)=\{y\in S_i(x)|\text{dist} (y, \partial_*S_i(x))\leq N^{11/16}\},\\
&\partial S_i^{(2)}(x)=\{y\in \Lambda\backslash S_i(x)|\text{dist} (y, \partial_*S_i(x))\leq N^{11/16}\}.\endaligned$$
We note that for $y\in\partial S_i^{(1)}(x)\cup\partial S_i^{(2)}(x)$, $\text{dist} (y, \partial_*S_i(x))\leq N^{11/16}\ll M$, the size of $\Lambda'$ cubes, and for $j$, $j'$ such that
$|j-j'|>N^{11/16}$, 
$$|W(j,j')|\ll e^{-N^{1/4}}.\tag 3.71$$

\item\item{$\bullet$} For all $x\in U$, with $\text{dist}(x,\partial_*S_{i-1})\geq M/4$, there exists $x'\in U$ such that 
$\Lambda'(x')\subset U$ and $\text{dist}(x, \partial_*\Lambda'(x'))\geq M/5$. We estimate $G_U(E)(x, y)$ with 
$x\in\partial S_{i-1}^{(2)}$, $y\in \partial S_{i+s}^{(1)}$ using $\Lambda'(x')\subset U$.

To obtain subexponential decay of off-diagonal elements $G_\Lambda(E)(x,y)$ we proceed as in \cite{BGS}.
This entails to estimate iteratively $G_{S_{n_i}}(E)(x, z)$, where $z\in\partial S_{n_i}^{(1)}$ and 
$A_{n_i}=S_{n_i}\backslash S_{n_i-1}$ is a good annulus. Let $n_i<m_i<n_{i+1}$, so that all the annulus 
in $S_{m_i}\backslash S_{n_i}$ are bad and all the annulus in  $S_{n_{i+1}}\backslash S_{m_i}=U$ are good.

Using the resolvent equation to relate $G_{S_{n_{i+1}}}(E)(x,z)$, $z\in\partial S_{n_{i+1}}^{(1)}$ with   
$G_{S_{n_{i}}}(E)(x,y)$, $y\in\partial S_{n_{i}}^{(1)}$, we have 
$$\aligned G_{S_{n_{i+1}}}&(E)(x,z)=\sum_{\Sb w\in U\\w'\in S_{n_{i+1}}\backslash U\endSb}G_{U}(E)(z,w)
\Gamma_U(w,w')G_{S_{n_{i+1}}}(E)(w',x)\\
&=\sum_{\Sb w\in \partial S_{m_{i}}^{(2)}\\w'\in \partial S_{m_{i}}^{(1)}\endSb}G_{U}(E)(z,w)\Gamma_U
(w,w')G_{S_{n_{i+1}}}(E)(w',x)+o(e^{-N^{1/4}})\\
&=\sum_{\Sb w\in \partial S_{m_{i}}^{(2)}\,w'\in\partial  S_{m_{i}}^{(1)}\\y\in \partial S_{n_{i}}^{(1)}\,y'\in \partial S_{n_{i}}^{(2)}\endSb}G_{U}(E)(z,w)\Gamma_U(w,w')G_{S_{n_{i}}}(E)(y,x)\Gamma_{S_{n_i}}(y,y')G_{S_{n_{i+1}}}(E)(y',w')\\
&\qquad\qquad\qquad\qquad+o(e^{-N^{1/4}}),
\endaligned \tag 3.72$$
where we used (3.71, 2.22). This is the analogue of (2.25, 2.26) in \cite{BGS}. Iterating (3.72) and taking the 
log log, lead to the conditions
$$\log\log[(e^{N^b})^{N/M}]<\log\log e^{N^{1/4}},\tag 3.73$$ 
$$\log\log e^{[M^{1/4}(N/M-N^b)]}>\log\log e^{N^{1/4}},\tag 3.74$$
for $N\gg 1$, where (3.73) originates from the estimates on $\Vert G_{S_{n_{i+1}}}\Vert$ and (3.74) from the decay estimates on $G_U$. 
(3.73, 3.74) in turn lead to $3/4+b<\tau<1-b$ for $M$, $N$ satisfying (3.69) (cf.\cite{BGS}).\hfill $\square$
\enddemo
In order to apply the above lemma, we need to convert the covering of $\Lambda=[-J^C, J^C]^{\nu+1}$ which is of 
diverging scales from $J^{\beta'}$ ($0<\beta'<1/5\alpha$) to $J^{C\beta}$ ($3/4<\beta<1$) to a covering of a single scale $J^{C\beta}$. 
Let $\Lambda'$ be cubes at scale $J^{C\beta}$ ($3/4<\beta<1$), i.e., 
$$\Lambda'(n,j)=(n,j)+[-J^{C\beta}, J^{C\beta}]^{\nu+1},\quad (n,j)\in\Bbb Z^{\nu+1}.\tag 3.75$$
Let $\Cal W$ be as in (3.62).
We define

\noindent$\Lambda'$ {\it bad} (for a fixed $\theta\in\Bbb R\backslash\Cal W(\Lambda)$), if
$\Lambda'\cap\tilde\Lambda\neq\emptyset$,  where $\tilde\Lambda$ is a bad $\Lambda_0$ box, i.e.,
(3.61) is violated.\hfill (3.76)

\noindent Let $\Gamma'$ be a covering of $\Lambda$ with $\Lambda'$ boxes defined in (3.75). Lemma 3.7
gives 
\proclaim{Lemma 3.9} Fix $\sigma\in(0, 1/5\alpha)$ as in Lemma 3.6.
For any fixed $\theta\in\Bbb R\backslash\Cal W(\Lambda)$, $\Gamma'$ has at most $\Cal
O(1)J^{5(\nu+1)}=(J^C)^{1-}\ll J^C$ ($5(\nu+1)<C\ll J^{\sigma/2}$) pairwise disjoint bad
$\Lambda'$ boxes. On the good $\Lambda'$ box, we have 
$$\align&\Vert G_{\Lambda'}(\theta)\Vert\leq e^{(J^{C\beta})^\sigma},\\
&|G_{\Lambda'}(\theta;n,j;n'j')|\leq e^{-|(n,j)-(n',j')|^{1/4}},\\
&\qquad\qquad \forall (n,j),\,(n'j')\,\text{such that }|(n,j)-(n',j')|>J^{C\beta}/10,\tag 3.77 \endalign$$
provided $$\max (5(\nu+1), \frac{\beta'}{{2\beta\sigma}})<C<\frac{\beta'}{\sigma}\ll J^{\sigma/2}.$$
So (3.61) is satisfied for the good $\Lambda'$.
\endproclaim 
\demo{Proof}
From (3.64), there are at most $\Cal O(1)J^{5(\nu+1)}$ points $(n,j)\in\Cal T$, such that 
$\Lambda_0(n,j)$ is a bad box, where the $\Cal O(1)$ only depends on the dimension $\nu$.
Since the $\Lambda'$ boxes are at scale $J^{C\beta}\gg J$, the scale of $\Lambda_0$ boxes,
there are at most $2^{\nu+1}$ non intersecting $\Lambda'$ boxes intersecting a given bad 
$\Lambda_0$ box. Hence for any fixed $\theta\in\Bbb R\backslash\Cal W(\Lambda)$, $\Gamma'$
has at most $\Cal O(1)J^{5(\nu+1)}$ pairwise disjoint bad $\Lambda'$ boxes with a slightly larger
$\Cal O(1)$.

On the good $\Lambda'$ boxes, we use the resolvent equation to obtain (3.77) as follows. We first
estimate $\Vert G_{\Lambda'}(\theta)\Vert$.
Let $\epsilon>0$ be arbitrary. We have
$$\aligned G_\Lambda(E+i\epsilon)(x,&y)=G_{U(x)}(E+i\epsilon)(x,y)\\
&+\sum_{\Sb z\in U(x)\\ z'\in \Lambda\backslash U(x)\endSb}G_{U(x)}(E+i\epsilon)(x,z)\Gamma_{U(x)}
(z,z')G_\Lambda(E+i\epsilon)(z',y),\endaligned$$
where $\Gamma_{U(x)}$ is as defined above (3.70), $U(x)$ is either $\Lambda_*$ or $\Lambda_0$ and $x\in U(x)$. Write $U$ for $U(x)$.
For $y'\in\partial_*U=\{y''\in U|\exists z\in\Lambda\backslash U,\, |z-y''|=1\}$, $x$ and $y'$ satisfy the distance 
condition in (3.42) if $U=\Lambda_*$ and the condition in (3.11) if $U=\Lambda_0$. It is easy to see 
that for all $x\in\Lambda$, such a $U(x)$ exists. 

Summing over $y$, we have 
$$\aligned \sum_y |G_\Lambda(E+i\epsilon)(&x,y)|\leq \sum_{y\in U(x)} \Vert G_{U}\Vert\\
&+2^{\nu+1}J^{C(\nu+1)}e^{-\frac{J^{\beta'}}{11}}\delta \sup_w \sum_y |G_\Lambda(E+i\epsilon)(w,y)|,\endaligned
$$
where we used (2.14, 2.22, 3.10, 3.11, 3.42, 3.43) and the conditions $\sigma,\,\beta'<1/5\alpha$, $\alpha>2\nu$.
Taking the supremum over $x$ and $U(x)$, we obtain the bound on $\Vert G_{\Lambda'}(\theta)\Vert$ in 
(3.77), provided $C>\beta'/2\beta\sigma$, $J\gg 1$.

Applying the resolvent equation one more time, similar to the proof of Lemma 3.8,  and using the bound on  $\Vert G_{\Lambda'}(\theta)\Vert$, we obtain the off-diagonal decay in (3.77). \hfill $\square$
\enddemo

Let $0<\beta'<1/5\alpha$, ($\alpha$ as in (3.23)), $3/4<\beta<1$. Choose $\sigma\in(0, 1/5\alpha)$
satisfying $5(\nu+1)<\beta'/\sigma$, which leads to 
$$0<\sigma<\frac{\beta'}{5(\nu+1)}<\frac{1}{25\alpha(\nu+1)}.$$
Choosing appropriate $C$ large enough so that both Lemme 3.8 and 3.9 are available, 
we then arrive at the following estimate for $G_\Lambda(\theta)=(H_\Lambda(\theta)-E)^{-1}$
valid for all $\Lambda=[-J,J]^{\nu+1}$, with $J$ large enough and any fixed $E$.
\proclaim{Proposition 3.10}
Assume $\omega$ is Diophantine, satisfying (3.23). For any $0<\beta'<1/{5\alpha}$, fix 
$$0<\sigma< \frac{\beta'}{20(\nu+1)}.$$ Then for all $0<\delta\ll 1$, there
exists $J_0$ such that the following statement is satisfied for all 
$\Lambda=[-J,J]^{\nu+1}$, with $J\geq J_0$: 
There exists $\Cal B(\Lambda, E)$ in
$\Bbb R$, with
$$\text{mes }\Cal B(\Lambda, E)\leq e^{-J^{\sigma/2}},\tag 3.78$$
such that if $\theta\in\Bbb R\backslash \Cal B(\Lambda, E)$, then 
$$\align&\Vert G_{\Lambda}(\theta, E)\Vert\leq e^{J^{\sigma}},\tag 3.79\\
&|G_{\Lambda}(\theta, E)(n,j;n'j')|\leq e^{-|(n,j)-(n',j')|^{1/4}},\\
&\qquad\qquad \text{for all } (n,j),\,(n'j')\,\text{such that }|(n,j)-(n',j')|>J/10.\tag 3.80\endalign$$
\endproclaim
\demo{Proof}
Choose $C>40(\nu+1)$ such that $$\frac{\beta'}{2\beta}<C\sigma<\beta'<1.$$
Assume $0<\delta\ll 1$. Let $J_0=|\log\delta|^{\frac{1}{C}}$.
For the scales $J_0$ to $J_0^C$, we use perturbation theory in $\delta$ as in Lemma 3.1,
with a further lowering of $\delta$ if necessary in order that (3.9-3.11) hold for all scales in 
$[J_0, J_0^C]$.
For the scales $\geq J_0^C$, Lemme 3.8 and 3.9 are available, in view of the choice of $C$.

Assume (3.78-3.80) hold at some scale $J\gg 1$. Using Lemma 3.9 in Lemma 3.8, we obtain the
corresponding estimate at scale $J^C$. Hence the proposition holds by induction.
\hfill$\square$
\enddemo
\bigskip
\head{\bf 4. Frequency estimates and the elimination of $E$.}\endhead
We now convert the estimate in $\theta$ in (3.78-3.80) for fixed $E$, $\omega$ into estimates in 
$\omega$ for 
fixed $\theta$ ($\theta=0$), eliminating the $E$ dependence in the process by excluding double 
resonances. This is the key estimate leading to the proof of the Theorem in section 5. 
We will first define a double resonant set in $(\omega, \theta)\in\Bbb T^\nu\times\Bbb R$.

We use $2$ scales $N$, $\bar N$ with $N<\bar N$, $\log\log\bar N\ll\log N$.
Let $\Lambda_N(j)=[-N,N]^{\nu+1}+(0,j)$, $j\in\Bbb Z$, $\Lambda_{\bar N}=[-\bar N,\bar N]^{\nu+1}$.
Assume $\Lambda_{\bar N}$ is resonant at $\theta=0$ for some $E$, i.e., $E$ is close to some eigenvalue of $H_{\Lambda_{\bar N}}$ (see (4.1) below). 
We prove in Lemma 4.1 that with  further reductions in the frequency set, the $\Lambda_N(j)$ boxes
with $|j|\leq\bar N$ and are at appropriate distances from 
$\Lambda_{\bar N}$ are all nonresonant  (at $\theta=0$). 

Toward that end, let $DC(\bar N)$ be the Diophantine condition to order $\bar N$, i.e., if $\omega\in DC(\bar N)$,
then $$\Vert n\cdot\omega\Vert_{\Bbb T}\geq\frac{c}{|n|^\alpha}\quad (\alpha>2\nu),\,\forall n\in[-\bar N, \bar N]^\nu\backslash\{0\},$$
where as before $\Vert\cdot\Vert_{\Bbb T}$ is the distance to the nearest integer.
Denote $H_{\Lambda_{\bar N}}$ by $H_{\bar N}$.  We define
$\Cal D(N, \bar N)\subset \Bbb T^\nu\times \Bbb R$ as
$$
\align
\Cal D(N, \bar N)&=\bigcup_{\{j\in\Bbb Z|\Lambda_N(j) \cap \Lambda_{\bar N}\not= \emptyset\}} 
\{(\omega, \theta)\in \Bbb T^\nu \times\Bbb R|\exists \, E, \text { such that }\\
&\Vert\big(H_{\bar N}(\omega, 0) -E\big) ^{-1} \Vert \geq e^{N^2} \text { and } \theta \in\Cal B
 \big(\Lambda_N(j), E\big)\},\tag 4.1
\endalign
$$
and for fixed $0<a<1$ and $C>1$, 
$$
S(N)=\bigcup_{\bar N\asymp  N^C} \big\{\Cal D(N, \bar N)
\bigcap (DC(\bar N)\times\Bbb R)\},\tag 4.2
$$
where $\bar N\asymp  N^C$ means $aN^C\leq\bar N\leq N^C/a$. We note that $a$ is kept fixed in 
this section.

The set $S(N)$ is the double resonant set in $(\omega, \theta)$ restricted to $\omega\in DC(\bar N)$,
$S(N)\subset \Bbb T^\nu\times \Bbb R$. It is the projection of the set $\tilde S(N)$, $\tilde S(N)\subset \Bbb T^\nu\times \Bbb R\times\Bbb R$, i.e., $S(N)=\text{Proj }_{\Bbb T^\nu\times \Bbb R}(\tilde S(N))$
with
$$
\tilde S(N)=\bigcup_{\bar N\asymp  N^C} \big\{\tilde \Cal D(N, \bar N)
\bigcap (DC(\bar N)\times\Bbb R\times\Bbb R)\},\tag 4.3$$
where
 $$
\align
\tilde \Cal D(N, \bar N)&=
\{(\omega, \theta, E)\in \Bbb T^\nu \times\Bbb R\times\Bbb R|
\Vert\big(H_{\bar N}(\omega, 0) -E\big) ^{-1} \Vert \geq e^{N^2}\\
&\quad\quad\text {and } \theta \in\Cal B
 \big(\Lambda_N(j), E\big) \text{ for some } \Lambda_N(j) \cap \Lambda_{\bar N}\not= \emptyset\},\tag 4.4
\endalign
$$
\proclaim {Lemma 4.1}
Let $N\in \Bbb N$ be sufficiently large, $C>1$ and $\delta>0$ sufficiently small. Let 
$\Omega_N$ be the union of $\Omega_{\bar N}^*$ over $\bar N\asymp  N^C$, where 
$\Omega_{\bar N}^*\subset\text{DC}(\bar N)$ is the set of $\omega$, such that 
\item\item{$\bullet$}
there exists $E\in\Bbb R$ such that
$$
\Vert\big(H_{\bar N}(\omega, 0)-E\big)^{-1} \Vert\geq e^{N^2};\tag 4.5
$$
\item\item{$\bullet$}
there exist $j\in \Bbb Z$, $\Lambda_N(j) \cap \Lambda_{\bar N}\not=\emptyset$, and $\ell\in \Bbb Z^\nu, |\ell|\asymp\exp[(\log N)^2]=K$
such that
$$
|(H_{\Lambda_N (j)}-E)^{-1}(\omega, \ell\cdot\omega)(m, m')|>e^{-|m-m'|^{1/4}}\tag 4.6
$$
for some $m, m' \in\Lambda_N(j), |m-m'|> N/10$. 

Then the set $\Omega_N$ satisfies
$$
\text {mes }\Omega_N\leq\exp[-\frac{1}{100}(\log N)^2] .\tag 4.7
$$
\endproclaim

The proof of Lemma 4.1 uses the following decomposition lemma \cite{B3, Lemma 9.9}.
\proclaim{Lemma 4.2}
Let $\Cal S\subset [0, 1]^{2\nu}$ be a semi-algebraic set of degree $B$ and $\text{\rm mes}_{2\nu} \Cal S
<\eta,
\log B\ll
\log 1/\eta$.
Denote by $(x, y)\in [0, 1]^\nu\times [0, 1]^\nu$ the product variable.
Fix $\epsilon>\eta^{1/2\nu}$.
Then there is a decomposition
$$
\Cal S =\Cal S_1 \bigcup \Cal S_2,
$$
with $\Cal S_1$ satisfying
$$
|\text{Proj}_x \Cal S_1|<B^C\epsilon\tag 4.8
$$
and $\Cal S_2$ satisfying the transversality property
$$
\text{\rm mes}_\nu(\Cal S_2\cap L)< B^C\epsilon^{-1} \eta^{1/2\nu},\tag 4.9
$$
for any $\nu$-dimensional hyperplane $L$ such that
$$
\max_{1\leq j\leq \nu}|\text{Proj}_L (e_j)|< \frac 1{100}\epsilon\tag 4.10
$$
where $e_j$ are the basis vectors for the $x$-coordinates.
\endproclaim
For our usage, the variable $(x,y)\in\Bbb R^{2\nu}$ in the lemma will be $(\omega,\theta)$ after
identifying $\theta$ with $(\theta,0)\in\Bbb R^{\nu}$.

\demo{Proof of Lemma 4.1}
Let  $(\omega, \theta)\in S(N)$, then $\theta \in\Cal B \big(\Lambda_N(j), E_k\big)$ for some $\Lambda_N(j)$ 
and some eigenvalue $E_k$ of $H_{\bar N}$ satisfying $|E_k-E|<e^{-N^2}$. This is because the perturbation
$e^{-N^2}$ essentially preserves the condition of the definition of a bad set, see (3.8-3.10) (c.f., proof of Lemma 4.1 \cite{BW1}). By the Diophantine restrictions on the frequencies 
$$\text{mes } S(N)\leq \Cal O(1)e^{-N^{\sigma/2}} N^{C(\nu+2)}\tag 4.11$$
from Proposition 3.10, where the third factor in the RHS is an upper bound on the number 
of $\Lambda_N(j)$ and the number of eigenvalues of $H_{\bar N}$.  

Since $\tilde S(N)$ is semi-algebraic with total degree at most $N^{C_1\nu}$ ($C_1>1$), cf., proof of Lemma 3.7, 
$S(N)$ its projection onto $\Bbb T^\nu\times\Bbb R$ is also semi-algebraic with degree at most 
$N^{C_2\nu}$ ($C_2>C_1$) \cite{Chap 9, B3}. (4.8) and Lemma 4.2 then conclude the proof by taking 
$\epsilon\asymp \exp[-(\log N)^2]$ (cf., also \cite{BW2}).

Here we also used the fact that the sum of the measure estimate in (4.9) over all $L$ corrsponding to the hyperplanes
$(\omega, \ell\cdot\omega)$ with $|\ell|\asymp\exp[(\log N)^2]$ is much smaller than $\exp[-(\log N)^2]$ for the 
corresponding $B$, $\eta$ and $\epsilon>0$.
\hfill$\square$
\enddemo

Lemma 4.1 shows that the boxes $\Lambda_N(\ell, j)=(\ell, j)+[-N, N]^{\nu+1}$, where $|\ell|\asymp  \exp[(\log N)^2]$, $j\in[-N^C, N^C]$ are non resonant with the box $\Lambda_{\bar N}=[-\bar N, \bar N]^{\nu+1}$,
$\bar N\asymp N^C$. Let $\Cal T$ be the tube region,
$$\Cal T=\{(n,j)\in\Bbb Z^{\nu+1}|j\in[-N^C, N^C]\}.\tag 4.12$$
The boxes $\Lambda_N(\ell, j)$ have centers in $\Cal T$. 

We now exclude resonances of $\Lambda_{\bar N}$ with boxes with centers in 
$\Bbb Z^d\backslash\Cal T$.
Recall from section 3, that boxes with centers in $\Bbb Z^d\backslash \Cal T$ are in fact rectangles:
$$\Lambda_*(n,2j)=(n,2j)+[-j^{\beta'},j^{\beta'}]^{\nu}\times[-j^\beta,j^\beta]\quad
(0<\beta'<1/5\alpha,\, 3/4<\beta<1),
\quad (n, 2j)\in\Bbb Z^d\backslash\Cal T.\tag 4.13$$
To exclude resonances of boxes $\Lambda_*(n,2j)$ where $\max (|n|, |2j|)\asymp \exp[(\log N)^2]$
with $\Lambda_{\bar N}$, we use direct perturbation. We have 
\proclaim{Lemma 4.3}
$$\aligned&\text{ mes }\{\omega\in[0, 2\pi)^\nu|\text { dist }(\sigma(H_{\bar N}),\sigma(H_{\Lambda_*}))\leq\kappa,\,\omega \text{ satisfies
(3.23)}\}\\
\leq&C\kappa|\Lambda_{\bar N}||\Lambda_*|,\qquad \qquad(\kappa\ll \exp[-\frac{1}{5}(\log N)^2])\endaligned\tag
4.14$$ 
for any $\Lambda_*=\Lambda_*(n, 2j)$ in (4.13) satisfying $|n|> \exp[\frac{1}{2}(\log N)^2]$
and $\max (|n|, |2j|)\asymp \exp[(\log N)^2]$, provided $0<\beta'<\min (1/20, 1/{5\alpha})$.
\endproclaim
\demo{Proof}
Let $\lambda_{m, k}$, $\phi_{m, k}$ be eigenvalues and eigenfunctions of $H_{\Lambda_*}$:
$$H_{\Lambda_*}\phi_{m, k}=\lambda_{m, k}\phi_{m, k}.$$
Write $\phi_{m, k}=\sum_{(m',k')\in\Lambda_*}a_{m',k'}\delta_{m',k'}$, with 
$\sum |a_{m',k'}|^2=1$. Then 
$$\lambda_{m,k}=\sum_{(m',k')}|a_{m',k'}|^2(m'\cdot\omega+k')+\delta\sum_{(m',k'),\,(m'', k'')}
a_{m',k'}\bar a_{m'',k''}\Delta_{m'm''}W_{k'k''}.$$

From (3.27) of Proposition 3.3 
$$\sum_{(m',k')\neq (m,k)}|a_{m',k'}|^2=\Cal O(\frac{1}{j^{1/20}}).\tag 4.15$$
The first order eigenvalue variation: 
$$\aligned m\cdot \frac{\partial}{\partial\omega}\lambda_{m,k}&=\sum_{(m',k')}|a_{m',k'}|^2 m\cdot m'\\
&=\sum_{(m',k')}|a_{m',k'}|^2 m\cdot m+\sum_{(m',k'),\, m'\neq m}|a_{m',k'}|^2 m\cdot (m'-m)\\
&\geq |m|^2-\frac{\Cal O(1)}{j^{1/20}}|m|\cdot j^{\beta'},\endaligned$$
where we used (4.15) and the fact that $|m-m'|\leq\Cal O(1)j^{\beta'}$ in $\Lambda_*$. Lowering
$\beta'$ to satisfy $\beta'<1/20$ if necessary, we obtain
$$m\cdot \frac{\partial}{\partial\omega}\lambda_{m,k}\geq \frac{1}{2}|m|^2\gg mN^C\tag 4.16$$
for $|m|>\exp[\frac{1}{2}(\log N)^2]\gg N^C$.

Let $\mu_{\ell,i}(\omega)$ be the eigenvalues of $H_{\bar N}$. Then $\mu_{\ell,i}(\omega)$ are piecewise
holomorphic in each $\omega_p$, $p=1,...,\nu$ and Lipshitz in $\omega$: $\Vert\mu_{\ell, i}\Vert_{\text{Lip}}
\leq C\bar N\asymp N^C$. Using this and (4.16),  we have that 
$$\aligned&\text{ mes }\{\omega\in[0, 2\pi)^\nu|\min_{\ell, i}\min_{m,k}|\lambda_{m,k}(\omega)-\mu_{\ell,i}(\omega)|
\leq\kappa\}\\
\leq&C\kappa|\Lambda_{\bar N}||\Lambda_*|.\endaligned$$
\hfill $\square$
\enddemo
Using Lemma 4.3, we arrive at 
\proclaim {Lemma 4.4}
Let $N\in \Bbb N$ be sufficiently large, $C>2/\beta'$,  $0<\beta'<\min (1/20, 1/{5\alpha})$ as in Lemma 4.3, $\delta$ sufficiently small. Let 
$\Omega'_N$ be a subset of the Diophantine set defined in (3.23) with the properties:

\item\item{$\bullet$}
there are $\bar N\asymp N^C $ and $E$ such that
$$
\Vert\big(H_{\bar N}(\omega, 0)-E\big)^{-1} \Vert\geq e^{N^2};\tag 4.17
$$
\item\item{$\bullet$}
there is  a $$\tilde\Lambda(n,2j)=(n,2j)+[-j^\beta, j^\beta]^{\nu+1}, \quad (3/4<\beta<1)$$
with $ |j|\geq N^C/2 $ and  $\max (|n|, |2j|)\asymp\exp[(\log N)^2] $ such that 
$$
|(H_{\tilde\Lambda(n,2j)}-E)^{-1}(\omega, 0)(m, m')|>e^{-|m-m'|^{1/4}}\tag 4.18
$$
for some $m, m' \in\tilde\Lambda$ satisfying $ |m-m'|> j^\beta/10$. 

Then
$$
\text {mes }\Omega'_N\leq e^{-N}.\tag 4.19
$$
\endproclaim

\demo{Proof}
We cover $\tilde\Lambda$ with $\Lambda_*$ of the form (4.13). We use a resolvent
expansion similar to the proof of Lemma 3.9 to obtain the opposite of (4.18). For that purpose we need to
take away an additional set of $\omega$ so that (3.42, 3.43) hold on all $\Lambda_*$ of the form
(4.13) at $\theta=0$. 

For any $\Lambda_*(n, 2j)$ such that
$|n|\geq \exp[\frac{1}{2}(\log N)^2] $,  $\max (|n|, |2j|)\asymp \exp[(\log N)^2]$, (4.14) is available.
Take one such $\Lambda_*$ and let $\kappa=e^{-2N}\gg e^{-N^2}$
for $N$ large. Assume (4.17) holds, so $|E-E_k|\leq e^{-N^2}\ll e^{-2N}$ for some
eigenvalue $E_k$ of $H_{\bar N}$. Lemma 4.3 then says that 
$$\Vert (E-H_{\Lambda_*})^{-1}\Vert\leq e^{2N}(1+2e^{-N^2})<2e^{2N}$$
by taking away a set in $\omega$ of measure $$\leq C\kappa|\Lambda_{\bar N}||\Lambda_*|\tag 4.20$$

Comparing (4.12) with (3.34) $J=N^C$ here. So  
$$\Vert (E-H_{\Lambda_*})^{-1}\Vert<2e^{2N}<e^{J^{\beta'/2}}=e^{N^{C\beta'/2}},$$
provided $C\beta'/2>1$ or $C>2/\beta'$. Therefore (3.42, 3.43) hold on this $\Lambda_*$.

Multiplying (4.20) by the number of all possible $\Lambda_*(n, 2j)$ such that  
$|n|\geq \exp[\frac{1}{2}(\log N)^2]$,  $\max (|n|, |2j|)\asymp \exp[(\log N)^2] $, we have that
(3.42, 3.43) are available on all such $\Lambda_*$ after taking away an additional set in $\omega$
of measure $\leq e^{-N}$.

For $\Lambda_*(n, 2j)$ such that $|n|<\exp[\frac{1}{2}(\log N)^2] $, since 
 $\max (|n|, |2j|)\asymp\exp[(\log N)^2] $,  $|j|\asymp\exp[(\log N)^2] $. So
 $$\text{dist }\big (E,\sigma(H_{\Lambda_*})\big)\asymp \exp[(\log N)^2], \quad \forall\omega\in[0,2\pi)^\nu,$$
 if $E$ satisfies (4.17). This is because $\Vert H_{\bar N}\Vert \leq\Cal O(\bar N)=N^C\ll \exp[(\log N)^2].$
 
 The proof now proceeds as the proof of Lemma 3.9 and we arrive at the conclusion. \hfill $\square$
\enddemo
\bigskip
\head{\bf 5. Proof of the Theorem}\endhead
We use Lemme 4.1, 4.4 to prove that the Floquet Hamiltonian $K$ in (1.15), or rather its
unitary equivalent $H$ in (1.19) has pure point spectrum. Let $\Omega_N$, $\Omega'_N$ 
be the frequency sets as in Lemme 4.1, 4.4, $N\in\Bbb N$, sufficiently large. Define
$\tilde\Omega$ to be
$$\tilde\Omega=\Bbb T^\nu\backslash\bigcup_{N\geq N_0}(\Omega_N\cup\Omega'_N), \quad N_0\gg1\text{ (depending on }\delta),\tag 5.1$$
then $\text{mes }\tilde\Omega\to (2\pi)^\nu$ as $\delta\to 0$.

Fix $\omega\in \tilde\Omega$. From the Schnol-Simon theorem \cite{CFKS, S}, to prove $H$ in (1.19) has pure point 
spectrum, it suffices to prove that the generalized eigenfunctions have fast decay, hence are in $\ell^2$.
More precisely, let $\psi$ be a non zero function on $\Bbb Z^{\nu+1}$ satisfying 
$$(H-E)\psi=0,\quad |\psi(m)|\leq 1+|m|^{c_0}\text{ for all }m\in\Bbb Z^{\nu+1},\tag 5.2$$ 
where $E$ is arbitrary and $c_0>0$ is some constant. We will prove using Lemme 4.1, 4.4 that
$\psi$ has subexponential decay and hence is in $\ell^2(\Bbb Z^{\nu+1})$.

We first verify that (4.5) is satisfied. This implies that (4.17) is also satisfied as they are the same condition. So we need to show that there is some box $\Lambda_{\bar N}$
centered at $0$, $\Lambda_{\bar N}=[-\bar N, \bar N]^{\nu+1}$, for some $\bar N\asymp N^C$
($C>2/\beta'$, $0<\beta'<\min (1/20, 1/{5\alpha})$ as in Lemme 4.1, 4.4), such that
$$
\Vert\big(H_{\bar N}-E\big)^{-1} \Vert\geq e^{N^2}.\tag 5.3
$$
For this we let $$T=\{(n,j)\in\Lambda_{\bar N}\,|\,j\in[-N^{C'}, N^{C'}]\},$$
$C'>1$ to be determined from (5.8-5.11), $C$ chosen to be $>C'$.
Let $\Lambda_{N^{C'}}$ be boxes of side length $2N^{C'}$.
We cover $T$ with $\Lambda_{N^{C'}}$, $\Lambda_{\bar N}\backslash T$ with $\Lambda_*$ of the form (4.13).

From Proposition 3.10, Lemma 3.7, there are at most ${\bar N}^{1-\tau_0}\quad (\tau_0>0)$ pairwise disjoint 
bad $\Lambda_{N^{C'}}$ boxes in $T$ by taking $C$ large enough. In $\Lambda_{\bar N}\backslash T$,
there is at most $1$ pairwise disjoint bad $\Lambda_*$ box by a further reduction in the $\tilde\Omega$ set as follows.

Assume $\exists (n, j)\in\Lambda_*$,  $\exists (m, k)\in\Lambda'_*$,  $\Lambda_*\cap\Lambda'_*=\emptyset$ such that 
$$
|\lambda_{n,j}(\omega)-\lambda'_{m,k}(\omega)|\leq\kappa,\tag 5.4$$
where $\lambda_{n,j}$,  $\lambda'_{m,k}$ are eigenvalues of $H_{\Lambda_*}$, $H_{\Lambda'_*}$.
From Proposition 3.3
$$\aligned \lambda_{n,j}&=n\cdot\omega+j+\Cal O(j^{-1/4})\\
\lambda'_{m,k}&=m\cdot\omega+k+\Cal O(k^{-1/4}).\endaligned$$
So $$|\lambda_{n,j}-\lambda'_{m,k}|\geq 1/2,$$
if $n=m$, since $(n,j)\neq (m,k)$. 

We only need to look at the case $n\neq m$.
This is similar to the proof of Lemma 4.3, except we look at the first order eigenvalue variation
in the ($n-m$) direction. Let 
$$\phi_{n,j}=\sum_{(n',j')\in\Lambda_*}a_{n',j'}\delta_{n',j'}$$ be the
eigenfunction with eigenvalue $\lambda_{n,j}$:
$$H_{\Lambda_*}\phi_{n,j}=\lambda_{n,j}\phi_{n,j};$$
and 
$$\psi_{m,k}=\sum_{(m',k')\in\Lambda'_*}b_{m',k'}\delta_{m',k'}$$ be the
eigenfunction with eigenvalue $\lambda'_{m,k}$:
$$H_{\Lambda'_*}\psi_{m,k}=\lambda_{m,k}\psi_{m,k}.$$
We have for the first order variation:
$$\aligned &(n-m)\cdot\frac{\partial}{\partial\omega}[\lambda_{n,j}-\lambda'_{m,k}]\\
=&(n-m)\cdot\frac{\partial}{\partial\omega}[\sum_{(n',j')\in\Lambda_*}|a_{n',j'}|^2n'\cdot\omega-\sum_{(m',k')\in\Lambda'_*}|b_{m',k'}|^2m'\cdot\omega]\\
=&(n-m)\cdot\frac{\partial}{\partial\omega}[(n-m)\cdot\omega+\sum_{(n',j')\in\Lambda_*,\, n\neq n'}|a_{n',j'}|^2(n'-n)\cdot\omega\\
&\qquad\qquad\qquad\qquad\qquad-\sum_{(m',k')\in\Lambda'_*,\, m\neq m'}|b_{m',k'}|^2(m'-m)\cdot\omega]\\
=&|n-m|^2+\Cal O(\frac{1}{j^{1/20}})\cdot j^{\beta'}|n-m|+\Cal O(\frac{1}{k^{1/20}})k^{\beta'}|n-m|\\
>&\frac{1}{2}|n-m|^2\\
\geq& \frac{1}{2},\endaligned\tag 5.5$$
where we used (3.27) of Proposition 3.3 and the definitions of $\Lambda_*$, $\Lambda'_*$ in (4.13).
So 
$$\text{mes }\{\omega||\lambda_{n,j}(\omega)-\lambda'_{m,k}(\omega)|\leq\kappa, \,\omega\text{ satisfies
}(3.23)\}\leq C\kappa.$$ 

Take $\kappa=e^{-(\log N)^2}$ and denote the set $\Omega''_N$ such that 
$\exists \Lambda_*\subset\Lambda_{\bar N}$, $\Lambda'_*\subset\Lambda_{\bar N}$, $\Lambda_*\cap\Lambda'_*=\emptyset$, $\exists (n, j)\in\Lambda_*$,  $\exists (m, k)\in\Lambda'_*$, such that (5.4) holds. Then 
$$\text {mes }\Omega''_N\leq e^{-\frac{1}{2}(\log N)^2},\tag 5.6$$ 
provided $2/\beta'<C\ll\log N$. 

Let $$\Omega{\overset\text{def }\to =}\tilde\Omega\backslash\bigcup_{N\geq N_0}\Omega''_N,\tag 5.7$$
$\text {mes }\Omega\to (2\pi)^{\nu}$ as $\delta\to 0$. We now assume $\omega\in\Omega$. Then in
$\Lambda_{\bar N}$, $\bar N\asymp N^C$, there are at most ${\bar N}^{1-\tau_0}$ ($0<\tau_0<1$)
bad $\Lambda_{N^{C'}}$ boxes with centers in $T$ and $1$ bad $\Lambda_*$ box with center in $\Lambda_{\bar N}\backslash T$.
Since $$\frac{\bar N}{N^{5(\nu+1)C'}}> {\bar N}^{\tau_0}\tag 5.8$$ for $C$ large enough, there has to be
an annulus $A$ at a distance $\asymp N^C$ to the origin of thickness $10 N^{C'}$ in $T$ and of thickness $N^C/4$ in $\Lambda_{\bar N}\backslash T$ (as there is at most $1$ bad $\Lambda_*$
box) devoid of bad points.

Since $(H-E)\psi=0$ from (5.2), let 
$$(H_{\tilde\Lambda(c)}-E)\psi=\xi_{\tilde\Lambda},\tag 5.9$$
where $\tilde\Lambda(c)$ is a box of type $\Lambda=\Lambda_{N^{C'}}$ or $\Lambda_*$ centered in $c$. In view of 
(1.20, 5.2), the restriction on $\beta$ for $\Lambda_*$ in (3.36)
$$|\xi_{\tilde\Lambda}(m)|\ll |m|^{c_0}e^{-c'[\text { dist }(m,\partial\tilde\Lambda)]^{1/2}},\tag 5.10$$
if $\text {dist}(m,\partial{\tilde\Lambda}^{(1)})>|m-c|_1$, where 
$$\partial{\tilde\Lambda}^{(1)}=\{y\in\tilde\Lambda|\exists z\in\Bbb Z^d\backslash\tilde\Lambda,\, |z-y|_1=1\}$$
and $|\,|_1$ is as defined above (3.46).

Let $c\in A$, $\tilde\Lambda(c)\subset A$. So $\tilde\Lambda(c)$ is a good box of type
$\Lambda_{N^{C'}}$ with centers in $T$ or $\Lambda_*$ of the form (4.13) with 
centers in $\Lambda\backslash T$ . Estimates (3.42, 3.80)
are available. From (5.9)
$$\psi=(H_{\tilde\Lambda(c)}-E)^{-1}\xi_{\tilde\Lambda}.$$
Using (5.10), (3.42) or (3.80), we have that 
$$|\psi(c)|<e^{-N^2},\tag 5.11$$
provided $C'>2/\beta'$ for all $c\in A$, such that $\exists \tilde\Lambda(c)\subset A$. We now choose
$2/\beta'<C'<C\ll \log N$, so that (5.8) is satisfied. (5.11) implies immediately that 
$$|\psi(m)|<e^{-N^2}\tag 5.12$$
for all $m\in A'\subset A$, where $A'$ is a smaller annulus contained in $A$ of thickness
$6N^{C'}$ in $T$ and thickness $N^C/5$ in $\Lambda_{\bar N}\backslash T$. 

Since $\text {dist } A'$ to the origin $\asymp N^C$, one can always find a square $\Lambda_{\bar N}$
centered at the origin, $\bar N\asymp N^C$, such that $\partial\Lambda_{\bar N}\subset A'$,
$$\align&\text { dist } (\partial\Lambda_{\bar N}, A'\cap T)>2N^{C'},\tag 5.13\\
&\text { dist } (\partial\Lambda_{\bar N}, A'\cap T^c)>N^{C}/10.\tag 5.14\endalign$$
Using (5.9) with $\Lambda_{\bar N}$ in place of $\tilde\Lambda(c)$, we have 
$$\psi=(H_{\bar N}-E)^{-1}\xi_{\bar N}.$$
Using (5.12-5.14), we have $\Vert\xi_{\bar N}\Vert<e^{-N^2}$ by slightly increasing
$C'$ if necessary.
Since we may always assume $\psi(0)=1$, (5.10, 5.11) give that 
$$\Vert (H_{\bar N}-E)^{-1}\Vert>e^{N^2}.$$

\noindent{\it Subexponential decay of eigenfunctions of} $H$.

Let $K=\exp[(\log N)^2]$. Lemme 4.1, 4.4 then imply (via a proof similar to the proofs of Lemme 3.8 and 3.9)
that the Green's function of the set
$$U{\overset\text{def }\to =}\Lambda_{2K}(0)\backslash \Lambda_{K}(0),$$
where $\Lambda_{L}(0)=[-L, L]^{\nu+1}$, exhibits off-diagonal decay, i.e., 
$$|G_U(E, m, n)|\leq e^{-|m-n|^{1/4}}\quad m,\,n\in U, \quad |m-n|>K/10.$$
Let $\xi_U$ be defined as in (5.9) with $U$ replacing $\tilde \Lambda$. 
Fix $m\in U$. Assume $\text{dist }(m,\partial U)\geq K/4$. For $n$ such that 
$\text {dist }(n,\partial U)>K/10$, we use the bound in (5.10) with $\xi_U$ replacing 
$\xi_{\tilde\Lambda}$. Otherwise we use the fact that $\xi_U$ is polynomially bounded
since $\psi$ is polynomially bounded. The equality
$$\psi(m)=\sum_n G_U(E, m,n)\xi(n),$$
then implies that 
$$|\psi(m)|\leq e^{-|m|^{\kappa}}\quad (0<\kappa<1/4),$$
provide $K$ and thus $|m|$ is large enough.
\hfill $\square$
\vfill\eject

\enddocument
\end